\RequirePackage{ifpdf}
\ifpdf 
\documentclass[pdftex]{sigma}
\else
\documentclass{sigma}
\fi

\usepackage[single]{accents}
\usepackage{mystyle}

\def\wb#1{\accentset{\hbox{\vrule height0.3pt width6pt depth0pt}}{#1}}

\usepackage{dsfont}
\def\1{{\mathds 1}}
\def\Z{{\mathds Z}}
\def\k{{\mathds k}}

\def\ol#1{\overline{#1}}

\numberwithin{equation}{section}
\numberwithin{theorem}{section}
\numberwithin{conjecture}{section}
\numberwithin{proposition}{section}
\numberwithin{lemma}{section}
\numberwithin{corollary}{section}
\numberwithin{definition}{section}
\numberwithin{example}{section}
\numberwithin{remark}{section}

\begin{document}

\allowdisplaybreaks

\renewcommand{\thefootnote}{$\star$}

\renewcommand{\PaperNumber}{057}

\FirstPageHeading

\ShortArticleName{On Griess Algebras}

\ArticleName{On Griess Algebras\footnote{This paper is a
contribution to the Special Issue on Kac--Moody Algebras and Applications. The
full collection is available at
\href{http://www.emis.de/journals/SIGMA/Kac-Moody_algebras.html}{http://www.emis.de/journals/SIGMA/Kac-Moody{\_}algebras.html}}}

\Author{Michael ROITMAN}
\AuthorNameForHeading{M. Roitman}

\Address{Department of Mathematics, Kansas State University, Manhattan, KS 66506 USA}
\Email{\href{mailto:misha.roitman@gmail.com}{misha.roitman@gmail.com}}

\ArticleDates{Received February 29, 2008, in f\/inal form July 28,
2008; Published online August 13, 2008}

\Abstract{In this paper we prove that for any commutative (but in
general non-associative) algebra $A$ with an invariant symmetric
non-degenerate bilinear form there is a graded vertex algebra
$V = V_0 \oplus V_2 \oplus V_3\oplus \cdots$, such that $\dim V_0 = 1$ and
$V_2$ contains $A$. We can choose $V$ so that if $A$ has a unit $e$, then $2e$ is the Virasoro element of $V$,
and if $G$ is a f\/inite group of automorphisms of $A$, then $G$ acts on $V$ as well.
In addition, the algebra $V$ can be chosen with a non-degenerate
invariant bilinear form, in which case it is simple.}


\Classification{17B69}

\renewcommand{\thefootnote}{\arabic{footnote}}
\setcounter{footnote}{0}

\section{Introduction}
A vertex algebra $V$ is a linear space, endowed with inf\/initely many bilinear products $(n):V\otimes V\to V$ and  a unit $\1\in V$,
satisfying certain axioms, see Section~\ref{sec:def}. In this paper we deal with
graded vertex algebras $V = \bigoplus_{i\in\Z}V_i$, so that $V_i (n) V_j
\subseteq V_{i+j-n-1}$ and $\1\in V_0$. A~vertex algebra is called OZ
(abbreviation of ``One-Zero'')  \cite{GNAVOA} if it is graded so that  $\dim V_0 =1$
and  $V_i = 0$  for $i=1$ or $i<0$. If $V$ is an OZ vertex
algebra, then \cite{flm} $V_2$ is a commutative (but not necessary associative)
algebra with respect to the product $(1):V_2\otimes V_2 \to V_2$, with an
invariant symmetric
bilinear form (i.e.\ such that $\form{ab}c = \form a{bc}$), given by the product
$(3):V_2\otimes V_2 \to V_0$. It is called {\em the
  Griess algebra of $V$}.\footnote{We note that the
term ``Griess algebra''  might not be the most successful one, as the
original Griess algebra~\cite{griess} is not quite a Griess algebra in
our sense.}

\subsection{Formulation of the results}
In this paper we prove the following result.
\renewcommand{\theenumi}{\alph{enumi}}

\begin{theorem}\label{main}  \ \null
\begin{enumerate}\itemsep=0pt
\item\label{main:emb}
For any commutative algebra $A$ with a symmetric
  invariant non-degenerate bilinear form there is a simple OZ vertex algebra $V$ such
  that $A\subseteq V_2$.
\item\label{main:vir}
If $A$ has a unit $e$, then $V$ can be chosen so that $\omega=2e$ is a
Virasoro element of $V$ (see Section~{\rm \ref{sec:def}} for the definition).
\item\label{main:aut}
If $G\subset \op{Aut}A$ is a finite group of automorphisms of
$A$, then  $V$ can be chosen so that $G\subset \op{Aut}V$.
\end{enumerate}
\end{theorem}
We prove this theorem under the assumption, that the ground f\/ield $\k$ is  a subf\/ield of
$\C$, since our proof uses some analytic methods (see Section~\ref{sec:regfun}). However, we believe that the
statement can be generalized to an arbitrary f\/ield of characteristic 0. Also, the
assumption that the form is non-degenerate does not seem to be very essential.

In fact we suggest that the following conjecture might be true:

\begin{conjecture} \ \null
\begin{enumerate}\itemsep=0pt
\item For any commutative algebra $A$ with a  symmetric
invariant bilinear form there is an OZ vertex algebra $V$ such that
$A=V_2$.
\item If $\dim A <\infty$, then $V$ can be chosen so that $\dim V_n
<\infty$ for $n=3,4,5,\ldots$.
\end{enumerate}
\end{conjecture}

It follows from Theorem~\ref{main}
that there are no Griess identities other than commutativity, in other words,
for any non-trivial identity in the variety of commutative algebras
with symmetric invariant bilinear forms there is a Griess algebra in which
this identity does not hold.

Here we outline our construction of $V$.
First we construct a vertex algebra $B=B_0\oplus B_2\oplus B_3 \oplus \cdots$, such that
$B_0$ is a polynomial algebra and $A\subset B_2$. In fact we  construct the vertex coalgebra
of correlation functions on $B$, def\/ined in Section~\ref{sec:coalgebras}, and then derive $B$ from it.
After that we f\/ind a suitable invariant bilinear form $\formdd$ on $B$ and set $V = B/\ker\formdd$.

We remark that our methods would perfectly work for a more general problem:
Given an ``initial segment'' $A_0\oplus A_1 \oplus \cdots \oplus A_m$ of a vertex algebra,
closed under those
of the vertex operations $(n)$ that make sense,  f\/ind a vertex algebra
$V=\bigoplus_{d\ge 0}V_d$ such that $V_d\supset A_d$ for $0\le d\le m$.

\subsection{Previously known results}
Probably the most famous example of OZ vertex algebras is the Moonshine module
$V^\natural$, constructed by Frenkel, Lepowsky and Meurman in
\cite{flm84,flm}, see also \cite{bor, bor92}. Its Griess algebra $V_2^\natural$ has
dimension 196\thinspace884, and dif\/fers from
the original 196\thinspace883-dimensional commutative algebra constructed by
Griess  \cite{griess} by having an additional identity element. The automorphism group of $V^\natural$ and of
$V_2^\natural$ is the Fischer--Griess Monster \cite{flm,griess,titsFG}.
It is proved by Dong et al.~\cite{dong_mon, dmzh},
that the representations of $V^\natural$ are completely reducible, and
the only irreducible representation is $V^\natural$ itself.
The present research was
primarily motivated by this construction.

Another example of OZ vertex algebra is a Virasoro vertex
algebra $\rm Vir$ \cite{fzh,wang}. It is generated by a single Virasoro
element $\omega \in {\rm Vir}_2$ so that the Griess algebra
of $\rm Vir$ is  $\k\,\omega$. The representation theory of the algebras $\rm Vir$
is investigated in \cite{wang}.

If $A$ is associative, than it is well known that $A$ can appear as a
Griess algebra, see \cite{lam96,xubook}. Lam \cite{lam99} also showed
the same for
a simple Jordan algebra of type $A$, $B$ or $C$. Other interesting
examples of OZ vertex algebras and their Griess algebras can be found
in \cite{GNAVOA}.

We remark that if a vertex algebra $V$ is graded so that
$V_n=0$ for $n<0$ and $\dim V_0 = 1$,
then $V_1$ is a Lie algebra with respect to the product $(0)$, with
invariant bilinear form given by product $(1)$.  The analogous problem
of f\/inding a vertex algebra $V$ such that $V_1$ is a given Lie algebra
has a well-know solution: to every Lie algebra $L$ with an invariant
bilinear form there corresponds an af\/f\/ine Lie algebra $\^L$, so that a
certain highest weight $\^L$-module has the desired vertex algebra
structure \cite{flm,fzh}.

\subsection{Organization of the manuscript}
In Section~\ref{sec:general} we recall some basic def\/initions and notations of the theory of vertex algebras.
Then in Section~\ref{sec:regfun} we consider a class of rational functions that we call {\em regular}. The
correlation functions of a suf\/f\/iciently nice vertex algebra will belong to that class.
Then in Section~\ref{sec:adm} we def\/ine a more narrow class of {\em admissible}
functions.
The correlation functions of the
algebras~$B$ and~$V$ that we construct later are admissible.
In Section~\ref{sec:coalgebras} we introduce a notion of vertex coalgebra of correlation functions, and show
how to reconstruct a vertex algebra by its coalgebra of correlation functions. In Section~\ref{sec:dg0}
we show that in some important cases the component of degree 0 of such vertex algebra is
isomorphic to a polynomial algebra.
In Section~\ref{sec:OZ} we study some easy properties of OZ vertex algebras, in particular (in Section~\ref{sec:vir})
investigate
the behavior of the correlation functions in the presence of a Virasoro element.   Then
in Section~\ref{sec:B} we construct certain vertex algebra $B$ using the coalgebra techniques developed in
Section~\ref{sec:coalgebras} and show (in Section~\ref{sec:proof}) how the existence of the algebra $B$ implies
Theorem~\ref{main}.

\subsection{Further questions}
Though the methods used in this paper are very explicit, it seems that the OZ vertex
algebras constructed here are of ``generic type'', i.e. they probably don't have these
nice properties people are looking for in vertex algebra theory~-- for example, an interesting
representation theory, various f\/initeness conditions, controllable Zhu algebra, etc.
It would be extremely interesting to recover the OZ vertex algebras mentioned above using
our approach, especially the Moonshine module $V^\natural$.

Also, it would be very
interesting to see whether any properties of the commutative algebra~$A$ (e.g.\ if $A$ is a
Jordan algebra) imply any properties of the OZ vertex algebra $V$, constructed in
Theorem~\ref{main}.

\section{General facts about vertex algebras}\label{sec:general}
Here we f\/ix the notations and give some minimal def\/initions. For more
details on vertex algebras the reader can refer to the books
\cite{flm,kac2,lepli,xubook}. Unless otherwise noted, we assume that all algebras and
spaces are over a ground f\/ield $\k\subset \C$.

\subsection[Definition of vertex algebras]{Def\/inition of vertex algebras}\label{sec:def}
\begin{definition}\label{dfn:vertex}
A vertex algebra is a linear space $V$ equipped with a family of
bilinear products $a\otimes b \mapsto a(n)b$, indexed by
integer parameter $n$, and with an element $\1\in V$, called the unit,
satisfying the identities (V1)--(V4) below.
Let $D:V\to V$ be the map def\/ined by $Da = a(-2)\1$.
Then the identities are: \smallskip
\begin{itemize}\itemsep=0pt
\item[\hypertarget{i}{(V1)}]
$a(n)b = 0$ for $n\gg 0$, \vskip7pt
\item[\hypertarget{ii}{(V2)}]
$\1(n)a = \delta_{n,-1}\,a$ \  \ and \ \
$a(n) \1 = \frac 1{(-n-1)!}\,D^{-n-1} a$,  \vskip8pt
\item[\hypertarget{iii}{(V3)}]
$D(a(n) b) = (Da)(n) b + a(n) (Db)$ \  \ and \ \
$(Da)(n) b  = - n\, a (n-1) b$,\vskip8pt
\item[\hypertarget{iv}{(V4)}]
$\displaystyle{a(m)\big(b(n) c\big) - b(n)\big(a(m) c\big) =
\smash[t]{\sum_{s\ge 0}\binom ms \big(a(s)b\big)(m+n-s)c}}$
\end{itemize}
for all $a,b,c\in V$ and $m,n\in \Z$.
\end{definition}

Another way of def\/ining vertex algebras is by using the generating series
\[
Y:V\to \op{Hom}(V, V\llp z\rrp)
\]
def\/ined for $a\in V$ by
\[
Y(a,z) = \sum_{n\in\Z}a(n)\,z^{-n-1},
\]
where $a(n):V\to V$ is the operator given by   $b\mapsto a(n)b$, and $z$ is a formal variable.
The most important property of these maps is that they are
{\em  local}: for any $a,b\in V$ there is $N\ge0$ such that
\begin{equation}\label{loc}
\big[ Y(a,w), Y(b,z)\big]\,(w-z)^N=0.
\end{equation}
In fact, this is the only essential condition that one needs to postulate to def\/ine vertex algeb\-ras~\cite{kac2,li}.
The minimal  number $N$ for which \eqref{loc} holds is called {\em
the locality} of $a$ and $b$, and is denoted by ${\rm loc}(a,b)$.
\begin{remark}
One could extend this def\/inition to allow a negative locality (see
\cite{cfva}), so that
\[
{\rm loc}(a,b) = \min\set{n\in\Z}{a(m)b=0\ \forall\; m\ge n}.
\]
\end{remark}

In terms of the series $Y$, the identities  (\hyperlink{ii}{V2}) and
(\hyperlink{iii}{V3}) read respectively
\begin{equation}\label{Y1}
Y(\1,z) = id,\qquad Y(a,z)\1 = \exp(Dz)a
\end{equation}
and
\begin{equation}\label{DY}
Y(Da,z) = \ad D{Y(a,z)} = \partial_z Y(a,z).
\end{equation}

Among other identities that hold in vertex algebras are  the quasi-symmetry
\begin{equation}\label{qs}
a(n)b = - \sum_{i\ge0} (-1)^{n+i} D^{(i)}
\big(b(n+i)a\big),
\end{equation}
and   the  associativity identity
\begin{gather}
\big(a(m)b\big)(n) c =
\sum_{s\ge0}(-1)^s \binom ms a(m-s)\big(b(n+s)c\big)\nonumber\\
\phantom{\big(a(m)b\big)(n) c =}{} -\sum_{s\le m}(-1)^s \binom m{m-s}
b(n+s)\big(a(m-s)c\big).\label{assoc}
\end{gather}

For $m\ge0$ this simplif\/ies to
\begin{equation*}
\big(a(m)b\big)(n) c =
\sum_{s=0}^m(-1)^s \binom ms \ad{a(m-s)}{b(n+s)}\,c,
\end{equation*}
which can also be derived from the identity (\hyperlink{iv}{V4})  of
Def\/inition~\ref{dfn:vertex} by some simple manipulations.

A vertex algebra $V$ is called {\em graded} (by the integers)
if $V = \bigoplus_{i\in\Z}V_i$ is a graded space, so that
$V_i(n)V_j\subseteq V_{i+j-n-1}$ and $\1\in V_0$.
It is often assumed that a vertex algebra $V$ is graded and
$V_2$ contains a special element $\omega$, called  the {\em Virasoro
  element} of $V$, such that $\omega(0)=D$,
$\omega(1)|_{V_i} = i$ and the coef\/f\/icients $\omega(n)$ generate a
representation of the Virasoro Lie algebra:
\begin{equation}\label{vir}
\ad{\omega(m)}{\omega(n)} = (m-n)\,\omega(m+n-1)+
\delta_{m+n,2}\,\frac 12\binom{m-1}3c
\end{equation}
for some constant $c\in \k$ called the {\em central charge} of $V$.
In this case $V$ is called {\em conformal vertex algebra} or, when $\dim V_i < \infty$, a {\em vertex
operator algebra}. The condition \eqref{vir} is equivalent to the following relations
\begin{equation}\label{virconf}
\omega(0)\omega = D\omega,\quad
\omega(1)\omega = 2\omega,\quad
\omega(2)\omega = 0,\quad
\omega(3)\omega = \frac c2,\quad
\omega(n)\omega = 0\quad\text{for}\ \ n\ge 4.
\end{equation}
This means that $\omega$ generates a {\em Virasoro conformal algebra}, see \cite{kac2}.

\begin{definition}[\cite{GNAVOA}]
A vertex algebra $V$ is called OZ (abbreviation of ``One-Zero'')
if it is graded so that  $V = \k \1 \oplus \bigoplus_{n\ge 2}V_n$.
\end{definition}
An OZ vertex algebra $V$ has $\dim V_0=1$ and $\dim V_1=0$, which
explains the name. The
component $V_2$ is a commutative (but not necessarily associative)
algebra with respect to the product $a\otimes b \mapsto ab = a(1)b$,
called the {\em Griess algebra} of $V$. The
commutativity follows from~\eqref{qs}. The algebra $V_2$ has a
bilinear form $\form ab = a(3)b$. From \eqref{qs} it follows that this
form is symmetric, and from \eqref{assoc} it follows that it it invariant:
$\form {ab}c = \form a{bc}$.


\begin{remark}
It should be noted that the idea that the 196\thinspace883-dimensional Griess's algebra
can be realized (after adjoining a unit) as a degree 2 component of a vertex algebra is due to Frenkel, Lepowsky
and Meurman \cite{flm}. The general fact that degree 2 component of any OZ vertex algebra
has a commutative algebra structure with a symmetric invariant bilinear form is mentioned
in this book as a triviality.
\end{remark}

One can def\/ine analogous structure on the components $V_0$ and $V_1$. Namely,
if a graded vertex algebra $V$ satisf\/ies $V_n=0$ for $n\le 0$, then $V_0$ is an
associative commutative algebra with respect to the product $(-1)$, and $V_1$ is a Lie
algebra with respect to the product $(0)$  with an invariant
symmetric bilinear form given by the product $(1)$.

We note that all def\/initions in this subsection make sense for $\k$ being a commutative
associative algebra containing $\Q$. In this case by ``linear space'' we understand a
torsion-free $\k$-module. This remark applies also to Sections~\ref{sec:uea}--\ref{sec:rad}
below, and to the parts of Section~\ref{sec:sl2} that does not refer to correlation functions.

\subsection{Correlation functions}\label{sec:corrfun}
Denote by $\Phi^l$ the space of rational functions in the variables
$z_1, \ldots, z_l$ of the form
\begin{equation}\label{PhiK}
\alpha=p(z_1,\ldots,z_l)
\prod_{1\le i<j\le l} (z_i-z_j)^{k_{ij}},\qquad
p\in \k[z_1,\ldots,z_l],
\end{equation}
where $k_{ij}\in \Z$.
Obviously we have a product $\Phi^l \otimes \Phi^m \to \Phi^{l+m}$ given by multiplying
the functions and renaming the variables.

Denote by $\ord_{ij}\alpha$ the order of $\alpha\in \Phi^l$ at $z_i-z_j$.
The space $\Phi^l=\bigoplus_{d\in\Z}\Phi^l_d$ is graded in
the usual sense, so that $\deg z_i = 1$.

Let $V$ be a graded vertex algebra, and let $f: V\to \k$ be a linear functional
of degree $d\in \Z$, i.e.\  $f(V_n)=0$ for $n\neq d$. Take some elements
$a_1,\ldots,a_l\in V$ of degrees $d_1,\ldots,d_l$ respectively and
formal variables $z_1,\ldots,z_l$. Consider the series
\begin{equation}\label{corrseries}
f\big(Y(a_1,z_1)\cdots
Y(a_l,z_l)\1\big)
= \sum_{m_1,\ldots,m_l\in\Z}
f\big(a_1(m_1)\cdots a_l(m_l)\1\big) \
z_1^{-m_1-1}\cdots z_l^{-m_l-1}.
\end{equation}

The following properties of the series \eqref{corrseries} can be deduced from  Def\/inition~\ref{dfn:vertex}
(see \cite{fhl}):
\begin{description}\itemsep=0pt
\item[\hypertarget{rat}{Rationality.}]
The series \eqref{corrseries} converge in the domain $|z_1|>\cdots > |z_l|$ to a rational function
$\alpha_f(z_1,\ldots,z_l)\in \Phi^l$ such that $\ord_{ij}\alpha_f\ge -
{\rm loc}(a_i,a_j)$. It is called a {\em correlation function}
of~$V$.
\item[\hypertarget{com}{Commutativity.}]
For any permutation $\sigma\in \Sigma_l$, the correlation function corresponding to \break
$a_{\sigma(1)}, \ldots, a_{\sigma(l)}$ and the same functional $f:V\to \k$ is
$\alpha_f(z_{\sigma(1)}, \ldots, z_{\sigma(l)})$.
\item[\hypertarget{assoc}{Associativity.}]
The series
\[
f\big(Y(Y(a_1,z_1-z_2)a_2, z_2)\,Y(a_3,z_3)\cdots Y(a_l,z_l)\1\big)
\]
converge in the domain $|z_2| > \cdots > |z_l|> |z_1-z_2|>0$ to $\alpha_f(z_1,\ldots,z_l)$.
\end{description}

Since $\deg a_i(m_i) =
d_i-m_i-1$, we get $\deg \alpha_f = \deg f -\sum_{i=1}^l d_i$.

It can be shown (see \cite{fhl,lepli}) that the rationality and commutativity properties
of correlation functions together with the conditions \eqref{Y1} and \eqref{DY}
can serve as an equivalent def\/inition  of vertex algebras.
We will use this fact in Section~\ref{sec:coalgebras} below.

In order to explain the meaning of the
\hyperlink{assoc}{associativity} condition, we need to introduce another def\/inition.
Take some $1\le i<j\le l$. A function $\alpha\in \Phi^l$ has expansion
\begin{equation}\label{alphaexpij}
\alpha(z_1,\ldots,z_n) = \sum_{k\ge k_0}
\alpha_k(z_1,\ldots,\^{z_i},\ldots,z_l)\, (z_i-z_j)^k
\end{equation}
for $\alpha_k\in \Phi^{l-1}$. Here and below the hat over a term indicates that this term is omitted. Then we def\/ine the
operators $\rho^{(k)}_{ij}:\Phi^l\to \Phi^{l-1}$ by setting
\begin{equation}\label{rho}
 \rho^{(k)}_{ij}\alpha = \alpha_k.
\end{equation}
An important property of these maps is that  for any $1\le i<j\le l$ and $1\le s<t\le l$,
such that $\{s,t\}\cap \{i,j\} = \varnothing$, and $m,k\in \Z$,
\begin{equation}\label{rhocom}
\rho^{(m)}_{st}\rho^{(k)}_{ij} = \rho^{(k)}_{ij}\rho^{(m)}_{st}.
\end{equation}

Now the \hyperlink{assoc}{associativity} condition means that if
$\alpha_f$ is a correlation function corresponding to the elements $a_1,\ldots,a_l$ and a
functional $f:V\to \k$, then
$\rho^{(k)}_{12}\alpha_f$ is the correlation function corresponding to the elements
$a_1(-k-1)a_2,\, a_3, \ldots, a_l$ and the same functional $f$.

\subsection[The action of $sl_2$]{The action of $\boldsymbol{sl_2}$}\label{sec:sl2}
In this paper  we will deal with vertex
algebras equipped a certain action of the Lie algebra $sl_2$.

\begin{definition}
A vertex algebra $V$ is said to have {\em $sl_2$ structure}, if
$V = \bigoplus_{d\in\Z}V_d$  is graded, and there is   a
locally nilpotent
operator $D^*:V\to V$  of degree $-1$, such that $D^*\1=0$ and
\begin{equation}  \label{adDst}
\ad{D^*}{a(m)} =  (2d-m-2)\,a(m+1) + (D^*a)(m)
\end{equation}
for every $a\in V_d$.
\end{definition}

Let $\delta: V\to V$ be the grading derivation,
def\/ined by $\delta|_{V_d} = d$. It is easy to compute that
if $D^*:V\to V$ satisf\/ies condition \eqref{adDst}, then
\[
\ad{D^*}D = 2\delta, \qquad
\ad{\delta}D = D,\qquad
\ad{\delta}{D^*} =- D^*,
\]
so that $D^*$, $D$ and $\delta$ span a copy of $sl_2$.

All vertex algebras in this paper are assumed to have
$sl_2$ structure, all ideals are stable under~$sl_2$ and homomorphisms of vertex algebras
preserve the action of $sl_2$ .

An element $a\in V$ such that $D^*a=0$ is called {\em minimal}.
It is easy to see that if  $V$ is generated by minimal elements, then
any operator $D^*:V\to V$ satisfying~\eqref{adDst} must be locally
nilpotent.

If $V$ has a Virasoro element $\omega$, then we can take $D^* =
\omega(2)$. Note that we always have
$D = \omega(0)$ and $\delta = \omega(1)$, therefore conformal vertex
algebras always have an~$sl_2$ structure.

Vertex algebras with an action of $sl_2$ as above were called
{\em quasi-vertex operator algebras} in~\cite{fhl} and homogeneous minimal
elements are sometimes called {\em quasi-primary}.

Now we describe the dual action on the correlation functions.
It follows from \eqref{DY} that the operator dual to $D$ is
$\Delta = \partial_{z_1} + \cdots + \partial_{z_l}$, so that
\[
f\big(DY(a_1,z_1)\cdots Y(a_l,z_l)\1\big) = \Delta\alpha_f(z_1,\ldots,z_l)
\]
for any homogeneous $a_1,\ldots,a_l\in V$ and $f:V\to \k$.
Note that $\Delta:\Phi^l \to \Phi^l$ is an operator of degree $-1$.

To describe the dual operator of $D^*$, consider the dif\/ferential operator
$\Delta^*(n,z) = z^2\partial_z +n\, z$. For the formal variables
$z_1,\ldots,z_l$, and for a sequence of integers $n_1, \ldots, n_l$ set
\[
\Delta^*(n_1,\ldots,n_l) = \Delta^*(n_1, z_1)+\cdots+\Delta^*(n_l, z_l).
\]
By \eqref{adDst}, we have
\[
f\big(D^*Y(a_1,z_1)\cdots Y(a_l,z_l)\1\big) = \Delta^*(2d_1,\ldots,2d_l)\, \alpha_f(z_1,\ldots,z_l)
\]
for minimal homogeneous elements $a_1,\ldots,a_l\in V$ of degrees $\deg a_i = d_i$
and a functional $f:V\to \k$.

Using the relations
\begin{gather}
\Delta^*(n_1,\ldots,n_l)\,(z_i-z_j)
= (z_i-z_j)\,\Delta^*(n_1,\ldots,n_i+1,\ldots,n_j+1,\ldots,n_l),\nonumber\\
\Delta^*(n_1,\ldots,n_l)\, z_i
= z_i\,\Delta^*(n_1,\ldots,n_i+1,\ldots,n_l),\label{Deltacom}
\end{gather}
where $z_i$ and $z_i-z_j$ are viewed as operators on $\Phi^l$, we see that $\Delta^*(n_1,\ldots,n_l)$ is an operator on~$\Phi^l$
of degree $1$.

We are going to need some easy facts about $sl_2$-module structure of $V$:
\begin{lemma}\label{lem:sl2} \ \null
\begin{enumerate}\itemsep=0pt
\item\label{sl2:neg}
If $d<0$, then  $V_d = (D^*)^{1-d}\, V_1$.
\item\label{sl2:DVm1}
$DV_{-1} \subseteq D^*V_1$.
\end{enumerate}
\end{lemma}
These statements hold for any graded $sl_2$-module $V$ on which $D^*$ is
locally nilpotent and $\delta\big|_{V_d}=d$ \cite{form}. The second statement follows easily from the f\/irst:
\[
DV_{-1} = DD^*V_0 = D^*DV_0 \subseteq D^*V_1.
\]

For vertex
algebras the action of $sl_2$ was also investigated in~\cite{dolima}.

\subsection{The universal enveloping algebra}\label{sec:uea}

For any vertex algebra $V$ we can construct a Lie algebra $L = \cff V$
in the following way \cite{bor,kac2,lihom,freecv}. Consider the linear space
$\k[t,t\inv]\otimes V$, where $t$ is a formal variable.
Denote $a(n) = a\otimes t^n$ for $n \in \Z$.
As a linear space,
$L$ is the quotient of $\k[t,t\inv]\otimes V$  by the subspace spanned
by  the  relations $(Da)(n) = -n\,a(n-1)$.  The brackets are
given by
\begin{equation}\label{com}
\ad{a(m)}{b(n)} = \sum_{i\ge0}\binom mi \big(a(i)b\big)(m+n-i),
\end{equation}
which is precisely the identity (\hyperlink{iv}{V4}) of Def\/inition~\ref{dfn:vertex}.
The spaces $L_\pm = \spn\set{a(n)}{n
\,\smash{\hbox to0pt{\raisebox{-4pt}{$<$}\hss}\raisebox{4pt}{$\ge$}}\,
0}\subset L$ are Lie subalgebras of $L$ so that $L = L_-\oplus L_+$.

\begin{remark}
The construction of $L$ makes use of only the products
$(n)$ for $n\ge 0$ and the map~$D$. This means that it works
for a more general algebraic structure, known as {\em conformal
algeb\-ra}~\cite{kac2,primc}.
\end{remark}

Now assume that the vertex algebra $V$ has an $sl_2$ structure.
Then (\hyperlink{iii}{V3}) of Def\/inition~\ref{dfn:vertex} and~\eqref{adDst} def\/ine derivations $D:L\to L$
and $D^*:L\to L$ so we get an action of $sl_2$ on $L$ by derivations.
Denote by $\^L = L\rtimes sl_2$ the corresponding semi-direct product.

The Lie algebra $\^L = \^{\cff V}$ and its universal enveloping algebra
$U=U(\^L)$ inherit the grading from $V$ so that
$\deg a(m) = \deg a -m-1$. The {\em Frenkel--Zhu topology} \cite{fzh}
on a homogeneous component $U_d$ is def\/ined by
setting the neighborhoods of 0 to be the spaces
$U^k_d = \sum_{i\le k} U_{d-i} U_i$, so that
\[
\cdots\subset U_d^{k-1} \subset U_d^k\subset U_d^{k+1}\subset \cdots \subset U_d,
\qquad \bigcap_{k\in \Z} U_d^k = 0,
\qquad \bigcup_{k\in \Z} U_d^k = U_d.
\]
Let
$\wb U = \bigoplus_{d\in\Z} \wb U_d$ be the completion of $U(\^L)$ in this
topology. Consider the ideal $I \subset \wb U$ generated by the relations
\begin{equation*}
\big(a(m)b\big)(n) =
\sum_{s\ge0}(-1)^s \binom ms a(m-s)b(n+s)
-\sum_{s\le m}(-1)^s \binom m{m-s}
b(n+s)a(m-s)
\end{equation*}
for all $a,b\in V$ and $m,n\in\Z$.
Note that the relations above are simply the associativity identity~\eqref{assoc}.
Denote by $W = \wb U/\bar I$ the quotient of $\wb U$ by the closure of
$I$.

For a f\/inite ordered set of elements
$\cal S = \{a_1,\ldots, a_l\}$, $a_i\in V$, let
$W_{\cal S}$ be the $\<D,\delta,D^*\>$-module generated by all
monomials
$a_1(m_1)\cdots a_l(m_l)\in W$, $m_i\in\Z$.

\begin{definition}[\cite{fzh,form}]\label{dfn:uea}
The universal enveloping algebra of $V$ is
\[
U(V) = \bigcup_{\cal S} \wb W_{\cal S}\subset W,
\]
where the union is taken over all f\/inite ordered sets
$\cal S\subset V$, and  $\wb W_{\cal S}\subset W$ is the completion of the space
$W_{\cal S}$ in the Frenkel--Zhu topology.
\end{definition}

\begin{remark}
In fact, it follows from  the  \hyperlink{com}{commutativity} property of correlation
functions (see Section~\ref{sec:corrfun})  that if $\cal S$ and $\cal S'$ dif\/fer by a permutation,
then $\wb W_{\cal S}= \wb W_{\cal S'}$.
\end{remark}

It is proved in~\cite{form} that any module over a vertex algebra $V$ is a continuous module over~$U(V)$,
in the sense that for any sequence $u_1,u_2,\ldots \in U(V)$
that converges to 0 and for any $v\in M$ we have $u_iv=0$ for $i\gg 0$.
Conversely, any $U(V)$-module $M$, such that $a(m)v=0$
for any $a\in V$, $v\in M$ and
$a(m)v=0$ for $m\gg0$, is a module over $V$.

\begin{remark}
The algebra $W = \wb U(\^L)/\bar I$ is also a good candidate for universal
enveloping algebra of~$V$. It has the following property \cite{fzh}:
consider a graded space $M$ such that
$M_d =0$ for $d\ll0$; then $M$ is a $W$-module if and only if $M$ is
an $V$-module.

On the other hand, we could def\/ine an algebra $\^U(V)$ such that any
series of elements from~$U(\^L)$, that make sense  as an operator on
any $V$-module, would converge in $\^U(V)$. However, this algebra is
too big for our purposes, for example there is no way of def\/ining an
involution in this algebra, as we do in Section~\ref{sec:forms} below.
\end{remark}

\subsection{Invariant bilinear forms}\label{sec:forms}
The key ingredient of our constructions is the notion of
invariant bilinear form on vertex algebra. Here we review the
results of~\cite{form}, that generalize the
results of  Frenkel, Huang and Lepowsky~\cite{fhl} and Li~\cite{liform}.

Let $V$ be a vertex algebra with an $sl_2$ structure, as in Section~\ref{sec:sl2}.
It is shown in~\cite{fhl,form}, that
there is an anti-involution $u\mapsto u^*$
on the universal enveloping algebra $U(V)$ such that $D\mapsto D^*$,
 $D^*\mapsto D$,
$\delta^*=\delta$ and
\begin{equation*}
 a(m)^* = (-1)^{\deg a} \sum_{i\ge 0} \frac 1{i!}\,
\big((D^*)^ia\big)(2\deg a-m-2-i)
\end{equation*}
for a homogeneous $a\in V$ and $m\in \Z$. In particular, if $D^*a=0$, then
\[
a(m)^*=(-1)^{\deg a}\,a(2\deg a-m-2),
\]
which can be written as
\begin{equation}\label{Yst}
Y(a,z)^* = \sum_{m\in \Z}a(m)^*\,z^{-m-1} = (-1)^{\deg a} \, Y(a,z\inv)\,z^{-2\deg a}.
\end{equation}

It is proved in \cite{form} that for any $u\in  U(V)_0$,
\begin{equation}\label{umust}
u\1-u^*\1\in D^*V_1.
\end{equation}

Let $K$ be a linear space over $\k$.
\begin{definition}[\cite{fhl,liform}]
A $K$-valued bilinear form $\form\cdot\cdot$ on $V$ is called {\it invariant} if
\[
\form {a(m)b}c = \form b{a(m)^*c}\qquad\text{and}\qquad
\form {Da}b = \form a{D^*b}
\]
for all $a,b,c\in V$ and $m\in\Z$.
\end{definition}
The radical
$\op{Rad}\form\cdot\cdot=\bigset{a\in V}{\form ab=0\ \forall\; b\in V}$
of an invariant form is an ideal of $V$. Also, since
$\form{\delta a}b =\form a{\delta b}$, we have $\form {V_i}{V_j}=0$
for $i\neq j$.

Given a $K$-valued invariant form $\form \cdot\cdot$ on $V$, one can consider a linear
functional $f: V_0\to K$ def\/ined by $f(a)= \form\1 a$. Since $f(D^*a)
= \form\1{D^*a} = \form {D\1}a = 0$, we get that
$f(D^*V_1)=0$. Also, the form can be reconstructed from $f$ by the
formula $\form ab = f\big(a(-1)^*b\big)$.

\begin{proposition}[\cite{liform,form}]\label{prop:forms}
There is a one-to-one correspondence between invariant $K$-valued
bilinear forms $\formdd$
on a vertex algebra $V$ and linear functionals $f:V_0/D^*V_1\to K$,
given by $f(a) = \form \1a$, $\form ab = f\big(a(-1)^*b\big)$.
Moreover, every invariant bilinear form on $V$ is symmetric.
\end{proposition}

\begin{remark}
We observe that a vertex algebra $V$ such that $V_0 = \k\1$ and
$D^*V_1=0$ is simple if and only if the invariant $\k$-valued bilinear form
on $V$ (which is unique by the above) is non-degenerate. Indeed, any
homomorphism $V\to U$ of vertex algebras must be an isometry, hence
its kernel must belong to the radical of the form.
\end{remark}

\subsection{Radical of a vertex algebra}\label{sec:rad}
Let $I=\<D^*V_1\>\subset V$ be the ideal of a vertex algebra $V$ generated
by the space $D^*V_1$.  Its degree 0 component $I_0 = U(V)_0D^*V_1$ is
spanned by the elements $a_1(m_1)\cdots a_l(m_l)D^*v$ such that
$a_i\in V$, $m_i\in\Z$, $\deg a_1(m_1)\cdots a_l(m_l) = 0$ and
$v\in V_1$, since we have $DV_{-1}\subset D^*V_1 \subset I_0$ by
Lemma~\ref{lem:sl2}\ref{sl2:DVm1}. Note that
Lemma~\ref{lem:sl2}\ref{sl2:neg} also implies that $V_d\subset I$ for $d<0$.

It follows from \eqref{qs} and \eqref{assoc} that $K=V_0/I_0$ is the
commutative associative algebra with respect to the product $(-1)$
with unit $\1$. Let $f:V_0 \to K$ be the canonical projection. By
Proposition~\ref{prop:forms}, the map $f$ corresponds to an invariant $K$-valued
bilinear form $\formdd$ on $V$.

\begin{definition}\label{dfn:rad}
 The {\em radical} of $V$ is $\rad V = \rad \formdd$.
\end{definition}

\begin{remark}
This def\/inition has nothing to do with the radical def\/ined in~\cite{dlmm}.
\end{remark}

Denote $\wb V = V/\rad V$. The following proposition summarizes some
properties of $\wb V$ that we will need later.

\begin{proposition}[\cite{form}]\label{prop:rad} \ \null
\begin{enumerate}\itemsep=0pt
  \item\label{rad:radrad}
$\rad(\wb V)=0$.
  \item\label{rad:cft}
$\wb V = \bigoplus_{n\ge 0} \wb V_n$, so that $\wb V_0 = V_0/I_0 = K$,
and $\wb V$ is a vertex algebra over $K$.
  \item\label{rad:ideals}
Every ideal $J_0 \subset K$ can be canonically extended to an ideal
$J\subset \wb V$, such that $J\cap \wb V_0 = J_0$.  The ideal $J$ is the maximal among
all ideals $I\subset \wb V$ with the property $I\cap \wb V_0 = J_0$.
In particular there are no non-trivial ideals $I \subset \wb V$ such
that $I\cap \wb V_0 = 0$.
  \end{enumerate}
\end{proposition}

The ideal $J\subset \wb V$ extending $J_0\subset \wb V_0$ is constructed  in the following
way: let $g:K\to K/J_0$ be the canonical projection, by Proposition~\ref{prop:forms}
it def\/ines a $K/J_0$-valued invariant bilinear form $\formdd_g$ on
$\wb V$. Then set $J = \rad\formdd_g$.

\section{Regular functions}\label{sec:regfun}
\subsection{Components}\label{sec:components}

Let $V$ be a vertex algebra with $sl_2$ structure. As in Section~\ref{sec:corrfun}, take some homogeneous elements
$a_1, \ldots, a_l\in V$ of
$\deg a_i = d_i$ and a functional $f:V_d \to \k$ of degree d, and let
$\alpha= \alpha_f(z_1, \ldots, z_l)\in \Phi^l$ be the corresponding correlation
function, given by~\eqref{corrseries}.  We have $\deg \alpha = d-\sum_i d_i$.

Denote by $\cal P$ the set of all partitions  $\{1,\ldots,l\}=I\sqcup J$  of the set
$\{1,\ldots,l\}$ into two disjoint subsets. For every $P = (I,J) \in \cal P$, the function
$\alpha$ has an expansion
\begin{equation}
  \label{comp}
  \alpha = \sum_{n\ge m} (\alpha)_n, \qquad\text{where}\qquad
(\alpha)_n = (\alpha)_n(P) = \sum_j \alpha'_{d-n,j} \, \alpha''_{n,j},
\end{equation}
for some $m\in \Z$. This expansion is obtained in the following way:
Let $I = \{i_1, \ldots, i_{|I|}\}$ and $J = \{j_1, \ldots,
j_{|I|}\}$. Expand $\alpha$ in power series in the domain
$|z_{i_1}|>|z_{i_2}| > \cdots > |z_{j_1}| > |z_{j_2}| > \cdots $ and
collect terms with powers of $\set{z_i}{i\in I}$ and $\set{z_j}{j\in J}$.
Note that the second sum in~\eqref{comp} is f\/inite. Here $\alpha'_{n,j}$ and  $\alpha''_{n,j}$ are rational functions depending on the variables
$\set{z_i}{i\in I}$, and  $\alpha''_{n,j}$ and $\set{z_i}{i\in J}$ respectively, and we have
\[
\deg\alpha'_{n,j} = n-\sum_{i\in I}d_i\qquad\text{and}\qquad
\deg \alpha''_{n,j} = n-\sum_{i\in J} d_i.
\]
We call the term  $(\alpha)_n$ in \eqref{comp} the  {\em component} of $\alpha$ of degree $n$
corresponding to partition $I\sqcup J$.
Note that $\alpha''_{n,j}\in \Phi^{|J|}$, while in general
$\alpha'_{n,j}\not\in \Phi^{|I|}$, since $\alpha'_{n,j}$ may have
a pole at $z_i$.


Assume that $\alpha= \alpha_f$  satisf\/ies $(\alpha)_n=0$ for $n<m$, and assume that
$I = \{i_1,\ldots, i_r\}$, $J = \{i_{r+1},\ldots,i_l\}$. Then $f\big(a_{i_1}(m_1)\cdots a_{i_l}(m_l)\1\big)=0$
whenever $\deg a_{i_{r+1}}\cdots a_{i_l}(m_l) < m$.
For example, suppose that $k$  is the order of $\alpha$ at $z_i =
z_j$. Take a partition
$\{1,\ldots,l\}=\{1,\ldots,\^i,\ldots,\^j,\ldots,l\}\sqcup \{i,j\}$.
Then $(\alpha)_n = 0$ for $n < d_i + d_j +k$, due to the \hyperlink{assoc}{associativity} property of Section~\ref{sec:corrfun}.

We are going to use the above terminology  even when $\alpha\in \Phi^l$ does not
necessarily correspond to a linear functional on a vertex algebra (for some f\/ixed integers
$k_1,\ldots, k_l$).

\subsection{Components of degree 0}\label{sec:comp0}
Denote by $\bar{\cal P}$ the set of {\em unordered} partitions of $\{1,\ldots,l\}$. Clearly, we
have a projection $\cal P\ni P\mapsto \bar P \in \bar{\cal P}$.

Fix some integers $d_1,\ldots, d_l$.
Suppose that a function $\alpha\in \Phi^l$ of degree $-\sum d_i$  satisf\/ies $(\alpha)_n(P) = 0$ for all $n<0$ and
$P\in \cal P$. Then the
expansion \eqref{comp} has a leading term $(\alpha)_0(P)$. It is easy to see that $(\alpha)_0$
depends only on the unordered partion $\bar P$.

\begin{proposition}\label{prop:comp0}
Suppose that for every partition $P=(I_1,I_2)\in \bar{\cal P}$ we have a function $\alpha(P) =
\sum_j \alpha^{(1)}_j\, \alpha^{(2)}_j$, where $\alpha_j^{(s)}$ depends on the variables
$\set{z_i}{i\in I_s}$, $\deg \alpha_j^{(s)} = -\sum_{i\in I_s}d_i$ and
$(\alpha^{(s)}_j)_d = 0$ for $d<0$, $s=1,2$. Assume that for any $Q\in \bar{\cal P}$ we have
\begin{equation}\label{alphaPQ}
(\alpha(P))_0(Q) = (\alpha(Q))_0(P).
\end{equation}
Then there is a function $\alpha\in
\Phi^l$, $k_P\in\k$, such that $(\alpha)_0(P) = \alpha(P)$.  Moreover, $\alpha$ is a
linear combination of $\alpha(P)$'s and their degree 0 components.
\end{proposition}

\begin{proof}
Introduce a linear ordering on the subsets of $\{1,\ldots,l\}$ such that $I<J$ if $|I|<|J|$,
and then extend it to $\bar{\cal P}$ so that $P=\{I_1,I_2\} < Q = \{J_1,J_2\}$ if $I_1< J_1$
and $I_1\le I_2$, $J_1\le J_2$. Set
$P_{\min} = \min\set{P\in \bar{\cal P}}{\alpha(P)\neq 0}$.
We will prove the existence of $\alpha$ by induction on $|\set{P\in\bar{\cal
    P}}{\alpha(P)\neq 0}|$. If $\alpha(P)=0$ for all $P$, take $\alpha=0$.

We observe that if $\alpha\in\Phi^l$ has degree $-\sum_i d_i$ and $(\alpha)_d = 0$ for
$d<0$, then the family of components $\{\alpha(P)=(\alpha)_0(P)\}_{P\in \bar{\cal P}}$ satisf\/ies~\eqref{alphaPQ}. Also, if collections $\{\alpha(P)\}$ and  $\{\beta(P)\}$ satisfy~\eqref{alphaPQ}, then so does  $\{\alpha(P) + \beta(P)\}$.

Now for any $P\in \bar{\cal P}_2$ set $\beta(P) = \alpha(P) - (\alpha(P_{\min}))_0(P)$. Obviously,
$\beta(P_{\min})=0$, and also, if $P{<}P_{\min}$, then $\alpha(P)=0$ and hence, using~\eqref{alphaPQ},
$\beta(P) = - (\alpha(P_{\min})_0)(P)= - (\alpha(P)_0)(P_{\min})=0$. By the above
observation, the collection $\{\beta(P)\}$ satisf\/ies~\eqref{alphaPQ}, therefore by induction, there is a
function $\beta \in \Phi^l$, such that
$(\beta)_0(P) = \beta(P)$ for any $P\in \bar{\cal P}$. Now take
$\alpha = \beta + \alpha(P_{\min})$.
\end{proof}

\begin{remark}
We can def\/ine components $(\alpha)_n(P)$ and the decomposition \eqref{comp} for partitions
$P$ of $\{1,\ldots,l\}$ into more than two parts.
Suppose we know the components $\alpha(P)$ for all such partitions $P$.
Then one can show that the function $\alpha\in \Phi^l$, such that $(\alpha)_0(P) = \alpha(P)$,
can be reconstructed by the following formula:
\begin{equation}\label{parts}
\alpha=  \sum_{P\in\bar{\cal P}}(-1)^{|P|}\,\big(|P|-1\big)!\,\alpha(P),
\end{equation}
where $\bar{\cal P}$ is the set of all unordered partitions of $\{1,\ldots,l\}$ and  $|P|$
is the number of parts of a~partition $P\in \bar{\cal P}$.
\end{remark}

\begin{remark}
It follows from the proof of Proposition~\ref{prop:comp0} that instead of \eqref{alphaPQ} it is enough to
require that the components $\alpha(P)$ satisfy the following property: If $\alpha(Q) = 0$
for some $Q\in \bar{\cal P}$, then $(\alpha(P))_0(Q)=0$ for every $P\in \bar{\cal P}$.
\end{remark}

\subsection{Regular functions}\label{sec:regular}
Recall that in Section~\ref{sec:sl2} we have def\/ined operators $\Delta$ and $\Delta^*$, so that
for a correlation function $\alpha_f(z_1,\ldots,z_n)\in \Phi^l$ corresponding to a
linear functional $f:V_d\to \k$ and elements $a_1\in V_{d_1},\ldots,a_l\in V_{d_l}$ we have
$f\big(DV_{d-1}\big)=0$ if and only if $\Delta \alpha_f = 0$ and
$f\big(D^*V_{d+1}\big)=0$
if and only if $\Delta^*(2d_1, \ldots, 2d_l)\,\alpha_f=0$.

It is easy to describe all functions $\alpha(z_1,\ldots,z_l)\in \Phi^l$, such that $\Delta \alpha = 0$.
These are the functions $\alpha$ that are invariant under translations, since{\samepage
\[
\alpha(z_1+t,\ldots,z_l+t) = \exp(t\Delta)\,\alpha(z_1,\ldots,z_l) = \alpha(z_1,\ldots,z_l).
\]
by the Taylor formula.
In other words, such $\alpha$ depends only on the dif\/ferences $z_i-z_j$.}

Now we will investigate the functions $\alpha$ which are killed by $\Delta^*$.

\begin{definition}\label{dfn:regular}
A function $\alpha\in \Phi^l$ is called {\em $(n_1, \ldots,  n_l)$-regular}
if  $\Delta^*(n_1, \ldots, n_l)\,\alpha_f=0$.
\end{definition}

\begin{example}\label{ex:S}
For an integer symmetric $l\times l$ matrix $\sf S = \{s_{ij}\}$ with
$s_{ii}=0$ def\/ine
\begin{equation}\label{pi}
\pi(\sf S) =
\prod_{1\le i<j\le l} (z_i-z_j)^{s_{ij}} \in \Phi^l.
\end{equation}

The relations \eqref{Deltacom} imply  that
\[
\Delta^*(n_1,\ldots,n_l)\big(\pi(\sf S)\big) =
\big((n_1+s_1)z_1 + \cdots + (n_l+s_l)z_1\big)\, \pi(\sf S),
\]
where   $s_i = \sum_j s_{ij}$.
Therefore, $\Delta^*(n_1,\ldots,n_l)\big(\pi(\sf S)\big) = 0$ if and only if
$s_i = -n_i$ for $i=1,\ldots,l$.  In this case the matrix $\sf S$ will
be called  $(n_1, \ldots, n_l)$-regular, so that  $\pi(\sf S)$ is a regular
function whenever $\sf S$ is a regular matrix.
\end{example}

\begin{remark}
One can show, though we will not use this here, that  the space of regular functions
$\ker \Delta^*(n_1,\ldots,n_l) \subset \Phi^l$ is spanned by
the products $\pi(\sf S)$ where $\sf S = \{s_{ij}\}$ runs over the set of $n_1,\ldots,n_l$-regular matrices such that $s_{ij} \ge k_{ij}$.
This description is analogous to the description of $\ker\Delta$ above.
Moreover, using this description, the dimensions of the homogeneous components of the spaces  $\ker
\Delta^*(n_1,\ldots,n_l) \subset \Phi^l$ can be given a combinatorial interpretation,
in fact, they  are certain generalizations of Catalan numbers.
\end{remark}

Assume we have a homogeneous linear functional $f:V\to \k$ such that $f(D^*V)=0$, and the
elements $a_i\in V_{d_i}$, $i=1,\ldots,l$. Then by Lemma~\ref{lem:sl2}\ref{sl2:neg}, we have $\deg f\ge
0$, and therefore the corresponding correlation function $\alpha_f$, given by
\eqref{corrseries}, is $(2d_1,\ldots,2d_l)$-regular and satisf\/ies $\deg \alpha_f \ge -\sum_i
d_i$. Moreover, if $\deg f = 0$, then also $f(DV)=0$ by Lemma~\ref{lem:sl2}\ref{sl2:DVm1}, and therefore
$\Delta\alpha_f = 0$.

Let us investigate the ef\/fect of the anti-involution  $u\mapsto u^*$ of the
enveloping algebra $U(V)$ (see  Section~\ref{sec:forms}) on the correlation functions.
Similarly to the series \eqref{corrseries}, one can consider the series
\[
f\big(\big(Y(a_1,z_1)\cdots Y(a_l,z_l)\big)^*\1\big) =
f\big(Y(a_l,z_l)^*\cdots Y(a_1,z_1)^*\1\big),
\]
which can be shown to converge in the domain $|z_1|<\cdots<|z_l|$ to a rational function\break
$\alpha_f^*(z_1,\ldots,z_l)\in \Phi^l$.  Since $f(D^*a_i) = 0$, we can
apply the formula \eqref{Yst} to each $a_i(n)^*$ and then it is easy ti
check that
\begin{equation}\label{alphastar}
\alpha_f^* = (-1)^{d_1+\cdots+d_l} \,
z_1^{-2d_1}\cdots z_l^{-2d_l}\alpha_f(z_1\inv,\ldots,z_l\inv).
\end{equation}
It follows from \eqref{umust} and the fact that $f(D^*V_1)=0$ that $\alpha^*_f = \alpha_f$.

In Section~\ref{sec:free}, given a collection of integers $\{k_{ij}\}$ for $1\le i<j\le l$,
we will construct a vertex algebra $F$ such that any function $\alpha \in
\Phi^l$ such that $\ord_{ij}\alpha\ge k_{ij}$ will be a correlation function on $F$, therefore the above properties hold
for any $(2d_1,\ldots,2d_l)$-regular function. It follows that the above properties of
correlation functions hold for all functions in $\Phi^l$.
Namely, we have the following proposition:

\begin{proposition}\label{prop:regular}
Let $\alpha \in \Phi^l$ be a $(2d_1,\ldots,2d_l)$-regular function. Then
$\deg \alpha \ge -\sum_i d_i$ and if
$\deg \alpha = -\sum_i d_i$, then $\Delta\alpha=0$ and $\alpha^*=\alpha$.
\end{proposition}
Here $\alpha^*$ is given by \eqref{alphastar}.
\begin{remark}
Proposition~\ref{prop:regular} can also be easily deduced from the fact that every regular function
is a linear combination of the products $\pi(\sf S)$. Also, one can show that a function
$\alpha\in \Phi^l$ of degree $-\sum d_i$ is $(2d_1,\ldots,2d_l)$-regular if and only if
$\alpha^*=\alpha$.
\end{remark}

\begin{corollary}\label{cor:regcomp}
Let $\alpha\in \Phi^l$ be a $(2d_1,\ldots,2d_l)$-regular function of degree $ -\sum_i
d_i$, and let $(\alpha)_n = \sum_j\alpha'_{-n,j}\, \alpha''_{n,j}$ be the degree $n$
component \eqref{comp} of $\alpha$ with respect to a partition
$\{1,\ldots,l\}=I\sqcup J$. Then the degree $n$ component of $\alpha$ with respect to
partition $J\sqcup I$ is $\sum_j(\alpha''_{n,j})^*(\alpha'_{-n,j})^*$.
\end{corollary}

\begin{remark}
It is easy to compute using \eqref{Yst} that for any correlation function $\alpha$ (and
therefore, for any function $\alpha\in \Phi^l$) one has
\[
\Delta^*(n_1,\ldots,n_l)\, \alpha^* = -(\Delta\alpha)^*,
\]
where
\[
\alpha^*(z_1,\ldots,z_l) = z_1^{-n_1}\cdots
z_l^{-n_l}\alpha(z_1\inv,\ldots,z_l\inv).
\]
As before, this can be easily computed without any reference to vertex algebras.
\end{remark}

\subsection{Admissible functions}\label{sec:adm}

In this section by ``regular'' we mean $(4,\ldots,4)$-regular, and set $\Delta^* =
\Delta^*(4,\ldots,4)$.

\begin{definition}\label{dfn:adm}
A regular function $\alpha\in \Phi^l$ is called {\em  admissible} if
for every partition $\{1,\ldots,l\} = I\sqcup J$ we have $(\alpha)_n=0$ for $n< 0$ or $n=1$.
If also  $(\alpha)_0 = 0$ for all partitions, then $\alpha$ is called {\em
indecomposable}.
\end{definition}
Denote the space of all admissible functions in $l$ variables by
$R^l\subset \Phi^l$, and the space of all indecomposable admissible functions
by $R_0^l\subset R^l$.

We have $R^1 = 0$,   $R^2 = \k (z_1-z_2)^{-4}$,  $R^3=R^3_0 = \k\, (z_1-z_2)^{-2} (z_1-z_3)^{-2} (z_2-z_3)^{-2}$, and it is easy to
compute, using e.g. the representation of regular functions by the products $\pi(\sf S)$,
that $\dim R^4_0 = 3$,  $\dim R^4 = 6$, $\dim R^5_0 = 16$,  $\dim R^5 = 26$ (compare with Section~\ref{sec:smalla} below).

We establish here a few simple properties of admissible functions. Recall that
the operators $\rho^{(k)}_{ij}:\Phi^l\to \Phi^{l-1}$ where def\/ined in \eqref{rho}.

\begin{proposition}\label{prop:adm}
Let $\alpha \in R^l$,  $l\ge3$ and $1\le i<j\le l$.
\begin{enumerate}\itemsep=0pt
\item \label{adm:ord}
$\ord_{ij}\alpha \ge -4$
and  if $\alpha\in R_0^l$, then $\ord_{ij}\alpha\ge -2$.
\item \label{adm:rho} $\rho^{(-4)}_{ij}\alpha\in R^{l-2}$,
$\rho^{(-3)}_{ij}\alpha=0$ and
$\rho_{ij}^{(-2)}\alpha \in R^{l-1}$.
\item \label{adm:poly}
The function $\alpha$ can be uniquely written as a linear combination of
the products of indecomposable admissible functions.
\end{enumerate}
\end{proposition}
The product in (\ref{adm:poly}) is understood in terms of the operation $\Phi^l\otimes
\Phi^m\to \Phi^{l+m}$ def\/ined in Section~\ref{sec:corrfun}.

\begin{proof}
To simplify notations, suppose $(i,j) = (l-1,l)$. Consider the
expansion \eqref{alphaexpij} for the function
$\alpha(z_1,\ldots,z_l)\in R^l$. If we expand every coef\/f\/icient
$\alpha_k(z_1,\ldots,z_{l-1})$ in the power series in $z_{l-1}$ around
$0$, we will get exactly the component expansion \eqref{comp} for the partition
$\{1,\ldots,l-2\}\sqcup\{l-1,l\}$. Then the minimal component is
\[
(\alpha)_{k_0+4} = \alpha_{k_0}(z_1,\ldots,z_{l-2},0)\, (z_{l-1}-z_l)^{k_0},
\]
where $k_0 =\ord_{l-1,l}\alpha$. This shows that $k_0\ge -4$ and,
since $(\alpha)_1=0$, we have
$k_0\neq -3$. Also, if~$\alpha$ is indecomposable, then $k_0\ge2$,
which proves  (\ref{adm:ord}).

Now assume that $k_0=-4$. Then we have
\[
0 = (\alpha)_1 = \Big(\, \frac{\partial \alpha_{-4}}{\partial z_{l-1}}\, \Big|_{z_{l-1}=0}\Big)\,
\big(z_{l-1}(z_{l-1}-z_l)^{-4}\big)
+ \big(\, \alpha_{_3}\big|_{z_{l-1}=0}\, \big) \, (z_{l-1}-z_l)^{-3}.
\]
Therefore, $\alpha_{-4}$ does not depend on $z_{l-1}$ and
$\alpha_3=0$ since $\alpha_3\in \Phi^{l-1}$. Since $\alpha$ does not
have components of negative degrees or of degree 1, neither do
$\alpha_{-2}$ and $\alpha_{-4}$. To prove  (\ref{adm:rho})
we are left to show that $\alpha_{-2}$ and $\alpha_{-4}$ are regular.

We have just seen that the
expansion \eqref{alphaexpij} for $\alpha$ has form
\begin{gather*}
  \alpha = \alpha_{-4}(z_1,\ldots,z_{l-2})\,(z_{l-1}-z_l)^{-4} +
\alpha_{-2}(z_1,\ldots,z_{l-1})\,(z_{l-1}-z_l)^{-2} + O\big((z_l-z_{l-1})^{-1}\big).
\end{gather*}
Applying $\Delta^*$ to this and using \eqref{Deltacom} we get
\[
0 = \Delta^*\alpha = \big(\Delta^*_2\alpha_{-4}\big)\,
(z_{l-1}-z_l)^{-4}
+\big(\Delta^*_1\alpha_{-2}\big)\,
(z_{l-1}-z_l)^{-2} + O\big((z_l-z_{l-1})^{-1}\big),
\]
where
$\Delta^*_s=\Delta^*(4,\ldots,4)$ ($l-s$ times), $s=1,2$,
which proves regularity of $\alpha_{-2}$ and $\alpha_{-4}$.

The proof of (\ref{adm:poly}) is very similar to the proof of Proposition~\ref{prop:comp0}.
Take a partition $P {=} (I_1,I_2){\in} \bar {\cal P}$.
We claim that if
$(\alpha)_0(P) = \sum_j \alpha'_{0,j} \, \alpha''_{0,j}$, then $\alpha'_{0,j}\in R^{|I_1|}$
and  $\alpha''_{0,j}\in R^{|I_2|}$. Indeed, we only need to check that $\alpha'_{0,j}$'s
and $\alpha''_{0,j}$'s are regular. Denote
$\Delta^*_s = \sum_{i\in I_s} \Delta^*(4,z_i)$, $s=1,2$ (see Section~\ref{sec:sl2}). Then
\[
0 = \Delta^*(\alpha)_0(P) = \sum_j (\Delta^*_1\alpha'_{0,j})\, \alpha''_{0,j} +
\sum_j \alpha'_{0,j}\, (\Delta^*_2\alpha''_{0,j}).
\]
Therefore, we see, using induction,  that $(\alpha)_0(P)$ is a linear combination of
products of indecomposable admissible functions. In particular we get $(\alpha)_0(P)\in R^l$.

Let $\alpha \in R^l$. If $(\alpha)_0(P)=0$ for every partition $P\in\cal P_2$, then
$\alpha\in R_0^l$. Otherwise, let $P\in \cal P_2$ be the minimal partition for which
$(\alpha)_0(P)\neq 0$. Then the
function $\beta = \alpha-(\alpha)_0(P)\in R^l$  will satisfy $(\beta)_0(Q)=0$ for all partitions $Q\le
P$. By induction, $\beta$ is a linear combination of  products of indecomposable
admissible functions, and hence so is $\alpha$.
\end{proof}

\begin{remark}
Alternatively, Proposition~\ref{prop:adm}\ref{adm:poly} follows from the formula \eqref{parts} in Section~\ref{sec:dg0} below.
\end{remark}

\begin{remark}
Suppose $\alpha\in \Phi^l$ is such that $\ord_{ij}\alpha=-2$. If $\alpha$ is regular, then so is
$\rho^{(-2)}_{ij}\alpha$. Indeed, applying $\Delta^*$ to
\[
\alpha = \sum_{k\ge -2}(z_i-z_j)^k\,  \rho^{(k)}_{ij}\alpha,
\]
we get, using \eqref{Deltacom},
\[
0 = \Delta^*\alpha = \sum_{k\ge -2} (z_i-z_j)^k\, \Delta^*(4, \ldots,  8+2k, \ldots,
4)\,\rho^{(k)}_{ij}\alpha.
\]
Here $8+2k$ stands at $i$-th position.
The coef\/f\/icient of $(z_i-z_j)^{-2}$ in the right-hand side is
$\Delta^*\rho^{(-2)}_{ij}\alpha$, which should be equal to 0.
In the same way one can check that if $\ord_{ij}\alpha=-4$, then $\rho_{ij}^{(-4)}\alpha$
is regular.
\end{remark}

Note that for $\alpha\in R^l$ and $1\le i<j\le l$ we have $\ord_{ij}\alpha \ge -4$ and
$\rho^{(-3)}_{ij}\alpha=0$.

\subsection{Admissible functions with prescribed poles}
In Section~\ref{sec:B} we will need the following property of admissible functions, which is reminiscent of the
Mittag--Lef\/f\/ler's theorem for analytic functions.

\begin{proposition}\label{prop:ML}
Let $l\ge 3$, and suppose that for each $1\le i<j\le l$ we fix admissible functions
$\alpha_{ij}^{(-2)}\in R^{(l-1)}$ and  $\alpha_{ij}^{(-4)}\in R^{(l-2)}$ satisfying the
following condition: For any $1\le s<t\le l$, such that $\{s,t\}\cap \{i,j\} = \varnothing$,
\begin{equation}\label{rhoalpha}
\rho^{(m)}_{st}\alpha^{(k)}_{ij} = \rho^{(k)}_{ij}\alpha^{(m)}_{st},\qquad
m,k=-2,-4.
\end{equation}
Then there exists a function $\alpha\in R^l$ such that $\rho^{(k)}_{ij}\alpha = \alpha^{(k)}_{ij}$
for all $1\le i<j\le l$ and $k=-2,-4$.
\end{proposition}
Note the similarity of the condition on $\alpha_{ij}^{(k)}$'s with \eqref{rhocom}.

In order to prove this proposition we need the following Lemma.

\begin{lemma}\label{lem:2poles}
Let $\alpha\in R^l$ be an admissible function.  Then for every $1\le i\le l$ one can
write $\alpha=\sum_m\alpha_m$ for some admissible functions $\alpha_m\in R^l$ that satisfy
the following properties:
\begin{itemize}\itemsep=0pt
\item[\rm (i)]
Either $\alpha_m = (z_i-z_j)^{-4}\beta$
for some $j\neq i$, where $\beta\in R^{l-2}$,
or $\ord_{ij}\alpha_m\ge -2$ for all $j\neq i$ and
\[
\big|\bigset{j\in \{1,\ldots,l\}\backslash\{i\}}{\ord_{ij}\alpha_m=-2}\big|\le2.
\]
\item[\rm (ii)]
For any $1\le s<t\le l$, if $\ord_{st}\alpha \ge -1$, then also $\ord_{st}\alpha_m \ge -1$,
and if $\ord_{st}\alpha =-2$, then $\ord_{st}\alpha_m \ge -2$.
\end{itemize}
\end{lemma}

In fact it will follow from the proof that if
$\ord_{ij}\alpha_m=\ord_{ik}\alpha_m=-2$ for some $j\neq k$, then $\alpha_m =
(z_i-z_j)^{-2}(z_i-z_k)^{-2}(z_j-z_k)^{-2}\beta$, where $\beta\in R^{l-3}$ does not depend
on $z_i, z_j, z_k$. Also, by Proposition~\ref{prop:adm}\ref{adm:poly} we can always assume that $\alpha_m$ is a product
of indecomposable admissible functions.

\begin{proof}
We use induction on $l$. If $l=2$ (respectively,
3), then $\alpha$ is a multiple of $(z_1-z_2)^{-4}$ (respectively,
$(z_1-z_2)^{-2}(z_1-z_3)^{-2}(z_2-z_3)^{-2}$) and we take $\alpha=\alpha_1$. So
assume that $l\ge 4$.

To simplify notations, assume that $i=1$. We also use induction on the number of multiple
poles of $\alpha$ as $z_1-z_j$,  $j=2,\ldots,l$, counting multiplicity.

Assume f\/irst that $\alpha$ has a pole of order 4 at one of $z_1-z_j$'s, which without loss of
generality we can assume to be $z_1-z_2$. Then set
\[
\gamma = (z_1-z_2)^{-4} \, \rho^{(-4)}_{12}\alpha.
\]
Obviously, $\gamma\in R^l$ and satisf\/ies (i). Since $\gamma$ does not have poles at
$z_1-z_j$, $z_2-z_j$ for $j\ge 3$, and~\eqref{rhocom} implies that $\ord_{ij}\gamma\ge
\ord_{ij}\alpha$ for all $3\le i<j\le l$,
$\gamma$  satisf\/ies (ii) as well.
Therefore, the function $\alpha'=\alpha-\gamma\in R^l$ has fewer
multiple poles at $z_1-z_j$. By induction, $\alpha' = \sum_m\alpha'_m$ for
$\alpha_m\in R^l$ satisfying (i) and (ii),  and hence $\alpha = \gamma + \sum_m\alpha'_m$.

Now assume that $\alpha$ has a double pole at some $z_1-z_j$, which is again can be taken
$z_1-z_2$. Then set $\beta(z_2,\ldots,z_l) = \rho^{(-2)}_{ij}\alpha\in R^{l-1}$.
By induction, we have $\beta = \sum_m \beta_m$, where the functions $\beta_m\in R^{l-1}$
satisfy conditions (i) and (ii) for $z_i=z_2$.
For each $\beta_m$, we need to consider two cases, that correspond to the dichotomy of the
condition (i):

{\bf Case 1.}
The function $\beta_m$ has a pole of order 4 at some $z_2-z_j$ for $j=3,\ldots,l$. Without loss of
generality, we can assume that $j=3$. Then $\beta_m = (z_2-z_3)^{(-4)}\beta'_m$ for some
$\beta'_m \in R^{l-3}$, and we set
\begin{equation}\label{gm1}
\gamma_m = (z_1-z_2)^{-2}(z_1-z_3)^{-2}(z_2-z_3)^{-2}\beta'_m.
\end{equation}

{\bf Case 2.}
The function $\beta_m$ has poles of orders at most 2 at all $z_2-z_j$ for
$j=3,\ldots,l$. Without loss of generality, we can assume that  $\ord_{2j}\beta_m\ge 1$
for $j\ge 5$. The we set
\begin{equation}\label{gm2}
\gamma_m = (z_1-z_2)^{-2}(z_1-z_3)\inv(z_1-z_4)\inv(z_2-z_3)(z_2-z_4)\beta_m.
\end{equation}

We need to show that in both cases the function $\gamma_m$ is admissible, satisf\/ies
conditions~(i) and~(ii) and $\rho^{(-2)}_{12}\gamma_m = \beta_m$. Indeed, assume that
these properties of $\gamma_m$ are established. Then set
$\alpha'=\alpha-\sum_m\gamma_m$. Since $\gamma_m$ satisf\/ies (ii), and
$\rho^{(-2)}_{12}\alpha' = 0$, the function $\alpha'$ will have less multiple poles in
$z_1-z_j$ than $\alpha$, therefore by induction, $\alpha'=\sum_m\alpha'_m$ for
$\alpha'_m\in R^l$ satisfying conditions~(i) and~(ii), and we take expansion $\alpha =
\sum_m \gamma_m + \sum_m\alpha'_m$.

Note that the conditions (i) and  $\rho^{(-2)}_{12}\gamma_m = \beta_m$ are obvious in both cases.

Case 1 is similar to the case when $\ord_{12}\alpha=-4$. We see that $\gamma_m$ is
admissible  by the def\/inition. Since $\ord_{23}\beta=-4$ and
$\ord_{12}\alpha=-2$, we must have $\ord_{13}\alpha\le -2$, which together with
\eqref{rhocom} establishes the property (ii) for $\gamma_m$.

So assume we are in Case 2. Condition (ii) follows from \eqref{rhocom} and the fact that the
only multiple pole of $\gamma_m$ at  $z_1-z_j$ and $z_2-z_j$ is at $z_1-z_2$.
We are left to show that $\gamma_m$ is admissible.

As it was mentioned above, we can assume that $\beta_m=\beta_{m1}\beta_{m2}\cdots$ is a
product of indecomposable admissible functions $\beta_{mt}\in R^{l_t}$, \
$\sum_tl_t=l-1$. Suppose $\beta_{l1}$ depends on $z_2$. Then $l_1>2$, since
$\ord_{23}\beta_m = -2$. Since our choice of $z_3$ and $z_4$ was based only on the
condition that $\ord_{2j}\beta_m \ge -1$ for $j\neq 3,4$, we can assume that $\beta_{m1}$
depends on $z_3$ and $z_4$. Therefore, in order to prove that $\gamma_m$ is admissible,
it is enough to show that if $\beta_m$ is indecomposable admissible,  then so is $\gamma_m$.

So assume that $\beta_m\in R_0^{l-1}$.
Applying $\Delta^*$ to $\gamma_m$ and using \eqref{Deltacom}
and the fact that $\beta_m$ is regular, we see that $\gamma_m$ is regular as well.
So we are left to verify that for every partition $\{1,\ldots,l\}= I\sqcup J$ we have
$(\gamma_m)_n=0$ for $n \le 1$.

Using Corollary~\ref{cor:regcomp}, we can assume without loss of generality that $1\in I$; otherwise we could
swap~$I$ and~$J$. Let
\[
(\beta_m)_n = \sum_j (\beta_m)'_{-n,j}\,(\beta_m)''_{n,j}
\]
be the component of $\beta_m$
corresponding to partition $\{2,\ldots,l\} = \big(I\backslash\{1\}\big)\sqcup J$. Since $\beta_m\in
R_0^{l-1}$, we have $(\beta_m)_n=0$ for $n\le 1$. We can expand the
factor
\[
\varkappa = (z_1-z_2)^{-2}(z_1-z_3)\inv(z_1-z_4)\inv(z_2-z_3)(z_2-z_4)
\]
in \eqref{gm2} as $\varkappa = \sum_s \varkappa'_s\, \varkappa''_s$,
where $\varkappa'_s$ depends on the variables $\bigset{z_i}{i\in I\cap \{1,2,3,4\}}$ and
$\varkappa''_s$ depends on the variables $\bigset{z_i}{i\in J\cap \{1,2,3,4\}}$, so that
$\deg \varkappa''_s\ge0$. We use here that $z_1$ appears in $\varkappa'_s$. Then the decomposition \eqref{comp} for
$\gamma_m$ becomes
\[
\sum_{n,j,s} \big(\varkappa'_s\,(\beta_m)'_{-n,j}\big)\, \big(\varkappa''_s\,(\beta_m)''_{n,j}\big),
\]
therefore, $(\gamma_m)_n=0$  for $n\le 1$.
\end{proof}

The proof of Proposition~\ref{prop:ML} is very similar to the proof of Lemma~\ref{lem:2poles}.

\begin{proof}[Proof of Proposition~\ref{prop:ML}]
We use induction on the number of non-zero functions among
$\set{\alpha_{ij}^{(k)}}{k=-2,-4,\, 1\le i<j\le l}$. If all of them are 0, then take
$\alpha =0$.

Assume f\/irst that some $\alpha_{ij}^{(-4)}\neq 0$. To simplify notations, we can take
$\alpha_{12}^{(-4)}\neq 0$. Then set
\[
\gamma = (z_1-z_2)^{-4} \, \alpha^{(-4)}_{12}\in R^l.
\]
As before, we see that $\gamma$ does not have poles at $(z_1-z_j), \, (z_2-z_j)$ for $j\ge
3$ and  $\ord_{ij}\gamma\ge
\ord_{ij}\alpha$ for  $3\le i<j\le l$, therefore the collection
$\{\alpha^{(k)}_{ij}-\rho^{(k)}_{ij}\gamma\}$ has fewer non-zero terms. This collection
satisf\/ies the condition \eqref{rhoalpha} because of the property \eqref{rhocom} of the
coef\/f\/icients $\rho_{ij}^{(k)}\gamma$. By induction,
there is a function $\alpha'$, such that $\rho^{(k)}_{ij}\alpha' =
\alpha^{(k)}_{ij}-\rho^{(k)}_{ij}\gamma$, and we can take $\alpha = \alpha'+\gamma$.

Now assume that $\alpha^{(-4)}_{ij}=0$, but $\beta = \alpha^{(-2)}_{ij}\neq0$, for some $1\le
i<j\le l$, which again can be assumed to be 1 and 2. Then
by Lemma~\ref{lem:2poles}, we can write $\beta(z_2,\ldots,z_l) = \sum_m\beta_m$ for some functions
$\beta_m\in R^{l-1}$ satisfying the conditions~(i) and~(ii) of Lemma~\ref{lem:2poles} for $z_i =
z_2$.
Exactly as in the proof of Lemma~\ref{lem:2poles}, without loss of generality we can consider two
cases for each $\beta_m$: when~$\beta_m$ has a pole of order 4 at $z_2-z_3$ and when~$\beta_m$ might have double poles at
$z_2-z_3$ and $z_2-z_4$ but at most simple poles at $z_2-z_j$ for $j\ge 5$. In each of
these cases def\/ine the function~$\gamma_m$ by the formulas~\eqref{gm1} and~\eqref{gm2}
respectively. As before, we see  each $\gamma_m$ is admissible, satisf\/ies the
property~(ii) of Lemma~\ref{lem:2poles} and $\rho_{12}^{(-2)}\gamma_m=\beta_m$. Therefore, setting
$\gamma=\sum_m\gamma_m$ as before, we see that the
collection  $\{\alpha^{(k)}_{ij}-\rho^{(k)}_{ij}\gamma\}$ has fewer non-zero terms and
satisf\/ies \eqref{rhoalpha}, so we f\/inish proof of the Proposition using induction as above.
\end{proof}

\section{The coalgebras of correlation functions}\label{sec:coalgebras}
\subsection{Spaces of correlation functions}\label{sec:algtocor}

Let $V=\bigoplus_d V_d$ be a vertex algebra with $sl_2$ structure. Assume that it has a set of
homogeneous generators $\cal G\subset V$ such that $D^*\cal G=0$.

\begin{remark}
The results in this section could be extended to the case when the
generators $\cal G$ are not necessarily minimal, but we do not need this
generalization here.
\end{remark}

Set $T(\cal G) = \set{a_1\otimes \cdots\otimes a_l\in V^{\otimes l}}{a_i\in \cal  G}$.
For any $\b a = a_1\otimes \cdots\otimes a_l \in T(\cal G)$, consider the space
\[
V^{\b a} = \spn_\k \bigset{a_1(n_1)\cdots a_l(n_l)\1}{n_i\in \Z}\subset V.
\]
Denote $V^{\b a}_d = V^{\b a}\cap V_d$.
The  \hyperlink{com}{commutativity} property of correlation
functions (see Section~\ref{sec:corrfun}) implies that
for any permutation $\sigma\in \Sigma_l$ and a  scalar $k\in \k$ we have
$V^{\sigma\b a} = V^{k \b a} = V^{\b a}$.

For a tensor $\b a = a_1\otimes \cdots \otimes a_l\in T(\cal G)$ and a subsequence
$I=\{i_1,i_2,\ldots\}\subset \{1,\ldots,l\}$ def\/ine
$\b a(I) = a_{i_1}\otimes a_{i_2}\otimes \cdots\in T(\cal G)$.

As it was explained in Section~\ref{sec:corrfun}, to any linear
functional $f:V^{\b a}_d\to \k$ we can correspond a~correlation function $\alpha_f \in
\Phi^l$ of degree $d-\sum_i \deg a_i$, such that $\ord_{ij}\alpha \ge -{\rm loc}(a_i,a_j)$.
Let
\[
\Omega^{\b a}=\bigoplus_d \Omega^{\b a}_d ,\qquad
\Omega^{\b a}_d = \set{\alpha_f}{f:V^{\b a}_d\to \k}\subset \Phi^l_{d-\sum\deg a_i}
\]
be the space of all such correlation functions, so that
$(V^{\b a}_d)^* \cong \Omega_d^{\b a}$.
\begin{definition}\label{dfn:vccf}
We will call the space
\begin{equation*}
\Omega = \Omega(V) = \bigoplus_{\b a \in T(\cal G)}\Omega^{\b a}
\end{equation*}
{\em the vertex coalgebra of correlation functions} of a vertex
algebra $V$.
\end{definition}
Note that $\Omega(V)$ depends on the choice of generators
$\cal G$, though we supress this dependence in the notation
$\Omega(V)$.
Also note that while each homogeneous component $\Omega^{\b a}\subset
\Phi^{|\b a|}$ consists of rational functions, the whole space
$\Omega(V)$ is not a subspace of $\Phi$.

The coalgebra structure on $\Omega(V)$, similar to the one def\/ined in \cite{hubbard}, is
manifested in the following properties, which easily follow from the properties of vertex
algebras (see Section~\ref{sec:general}):

\renewcommand{\theenumi}{$\Omega$\arabic{enumi}}
\begin{enumerate}\itemsep=0pt
\setcounter{enumi}{-1}
\item\label{O:1}
$\Omega^1 = \k$, $\Omega^a = \k[z]$ for every
$a\in \cal G$,
and for $\b a = a_1\otimes\cdots\otimes a_l\in T(\cal G)$,  $l\ge 2$.
\item\label{O:finite}
$\ord_{ij}\alpha \ge -{\rm loc}(a_i,a_j)$ for any $\alpha\in \Omega^{\b a}$.
\item\label{O:sym}
$\Omega^{\b a} = \sigma\Omega^{\sigma\b a}$ for any permutation $\sigma\in \Sigma_l$.
\item\label{O:Delta}
The space $\Omega^{\b a}$ is closed under the operators
$\Delta = \sum_i\partial_{z_i}$ and $\Delta^* = \sum_i (z_i^2\,\partial_{z_i} + 2(\deg a_i)\, z_i)$.
\item\label{O:corr}
Set $\b b = a_2\otimes\cdots\otimes a_l\in T(\cal G)$. Then
any function  $\alpha\in \Omega^{\b a}$ can be expanded at $z_1=\infty$ into a series
\begin{equation}\label{iter}
\alpha(z_1,\ldots, z_l)=
\sum_{n\ge n_0} z_1^{-n-1} \alpha_n(z_2,\ldots z_l),
\end{equation}
where $\alpha_n \in \Omega^{\b b}$.
\end{enumerate}

The action of $\Sigma_l$ on $\Omega^{\b a}$ in (\ref{O:sym}) is def\/ined by
$(\sigma \alpha)(z_1.\ldots,z_l) =
\alpha(z_{\sigma(1)},\ldots,z_{\sigma(l)})$, so that (\ref{O:sym}) is
just the  \hyperlink{com}{commutativity} property of Section~\ref{sec:corrfun}. It
implies that the space $\Omega^{\b a}$ for $\b a =
a_1\otimes\cdots\otimes a_l\in T(\cal G)$ is symmetric under the group $\Gamma_{\b
  a}\subset \Sigma_l$ generated by all transpositions $(i\ j)$ whenever $a_i=a_j$.

Note that in order to get the expansion \eqref{comp} of a function $\alpha\in \Omega^{\b a}$,
we need to apply a suitable permutation to the variables $z_1,\ldots,z_n$, and then iterate
the expansion \eqref{iter} several times. Combining this observation with the property
(\ref{O:sym}), we see that (\ref{O:corr}) can be reformulated as follows:
\renewcommand{\theenumi}{$\Omega$\arabic{enumi}$'$}
\begin{enumerate}
\setcounter{enumi}{3}
\item\label{O:corr1} For a partition $\{1,\ldots,l\}=I\sqcup J$, denote $\b a'' = \b
a(J)$. Then the component of degree $n$ of a~function $\alpha(z_1,\ldots,z_l) \in
\Omega^{\b a}$ of degree $d-\sum_i\deg a_i$ can be written as  $(\alpha)_s =
 \sum_j\alpha'_{d-n,j}\alpha''_{n,j}$ so that
 $\alpha''_{n,j}\in \Omega^{\b a''}$.
\end{enumerate}
\renewcommand{\theenumi}{\alph{enumi}}

\subsection{Universal vertex algebras}\label{sec:cortoalg}
Now we want to present a converse construction: given a space of functions $\Omega$,
satisfying the conditions (\ref{O:1})--(\ref{O:corr}), we will construct a vertex
algebra $V = V(\Omega)$, such that $\Omega = \Omega(V)$.

 Let $\cal G$ be a set. For any $a\in \cal
G$ f\/ix its degree $\deg a\in \Z$, and  for any pair
$a,b\in \cal G$ f\/ix a number ${\rm loc}(a,b)\in \Z$.

\begin{theorem}\label{thm:V}
Let $\Omega = \bigoplus_{\b a\in T(\cal G)}
\Omega^{\b a}$ be a graded space constructed from rational functions
as above, satisfying conditions
{\rm (\ref{O:1})--(\ref{O:corr})}. Then there exists a vertex algebra
\[
V =V(\Omega)=  \bigoplus_{\lambda \in \Z_+[\cal G]}V^\lambda,
\]
generated by $\cal G$ so that $a\in V_{\deg a}$, $D^*a=0$ for any $a\in \cal G$,
such that $\Omega = \Omega(V)$ is the vertex coalgebra of correlation
functions of $V$ (see Definition~{\rm \ref{dfn:vccf}}).
\end{theorem}

\begin{proof}
Let $\Omega^{\b a}_d\subset \Omega^{\b a}$ be the subspace of functions of degree
$d-\sum_i\deg a_i$. The condition (\ref{O:finite}) implies that
$\Omega^{\b a}=\bigoplus_d \Omega^{\b a}_d$  so that $\Omega^{\b a}_d = 0$ when $d\ll 0$
and $\dim \Omega^{\b a}_d< \infty$.

For each $\b a \in T(\cal G)$ set $V^{\b a} = (\Omega^{\b a})'$ to be the graded dual
space of~$\Omega^{\b a}$. We def\/ine degree on~$V^{\b a}$ by setting $\deg v = d + \sum_i \deg a_i$ for $v: \Omega^{\b a}_d \to \k$.
For a permutation $\sigma
\in \Sigma_l$ we identify~$V^{\b a}$ with $V^{\sigma \b a}$ using (\ref{O:sym}). In this
way for every $\lambda =a_1+\cdots + a_l\in \Z_+[\cal G]$ we obtain a space $V^\lambda =
V^{a_1\otimes\cdots\otimes a_l}$, and set $V = V(\Omega) = \bigoplus_{\lambda \in \Z_+[\cal G]}V^\lambda$.

For every $\lambda = a_1 + \cdots + a_l\in \Z_+[\cal G]$
choose a basis $\cal B^\lambda_d$ of $V^\lambda_d$. Set $\cal B^\lambda = \bigcup_d
\cal B^\lambda_d$. Let $\set{\alpha_u}{u\in \cal B_d^\lambda}$ be the dual basis of
$\Omega^{\b a}_d$, where $\b a = a_1 \otimes\cdots\otimes a_l\in T(\cal G)$, so that
$\deg \alpha_u = \deg u -\sum_i \deg a_i$.
For a permutation $\sigma\in \Sigma_l$ the set $\set{\sigma\alpha_u}{u\in \cal
  B_d^\lambda}$ is the basis of $\Omega_d^{\sigma \b a}$, dual to $\cal B^\lambda_d$.

We choose these bases so that $\cal B^0 = \{\1\}$ and $\alpha_\1 = 1\in \Omega^0 =
\k$. Also, for a generator $a\in \cal G$ of degree $d$ we have $\dim V^a_d = 1$, since
$V^a_d  = (\Omega^a_0)^*$ and $\Omega^a_0=\k$ due to (\ref{O:1}). We can identify the
only element of $\cal B^a_d$ with $a$ so that $\alpha_a = 1\in \Omega^a_0$.

We def\/ine the operators $D:V^\lambda\to V^\lambda$ and $D^*:V^\lambda\to V^\lambda$ as the
duals to
$\Delta:\Omega^{\b a}\to \Omega^{\b a}$ and
$\Delta^*(2\deg a_1,\ldots, 2\deg a_l):\Omega^{\b a}\to \Omega^{\b a}$
respectively (see Section~\ref{sec:sl2}). Since $\Delta^*(2\deg a) \Omega^a
\subset z \k[z]$, we have $D^*a=0$ for every $a\in \cal G$.

Now we are going to def\/ine vertex algebra structure  $Y:V\to \op{Hom}(V,V\llp z\rrp )$ so
that for any $\lambda = a_1 + \cdots + a_l\in \Z_+[\cal G]$ we have
\begin{gather}\label{longY}
Y(a_1,z_1)\cdots Y(a_l,z_l)\1 = \sum_{u\in \cal B^\lambda}\alpha_u(z_1,\ldots,z_l)\,u,\\
\label{Yassoc}
Y(a_1,z_1+z)\cdots Y(a_l,z_l+z)\,w
= Y\big(Y(a_1,z_1)\cdots Y(a_l,z_l)\1, z\big)\,w.
\end{gather}
These identities are to be understood in the following sense. The left-hand side of~\eqref{longY} converges to the $V$-valued rational function on the right-hand side in the
region $|z_1|> |z_2| > \cdots > |z_l|$. The left and right-hand sides of~\eqref{Yassoc}
converge to the same $V$-valued rational function in the regions
$|z_1+z|> |z_2+z| > \cdots > |z_l+z|$
and $|z| > |z_1|> |z_2| > \cdots > |z_l|$
respectively.

Take, as before,  $\lambda = a_1 + \cdots + a_l\in \Z_+[\cal G]$ and
$\b a = a_1 \otimes\cdots\otimes a_l\in T(\cal G)$. Let $a\in \cal G$ be a~generator  of
degree $d$.
First we def\/ine the action of $Y(a,z)$ on $V^\lambda$.

For any $v\in \cal B^{a+\lambda}$
expand the corresponding basic function $\alpha_v$ as in \eqref{iter}:
\[
\alpha_v(z,z_1, \ldots,z_l) = \sum_n z^{-n-1}\, \alpha_n(z_1,\ldots,z_l),\qquad
\alpha_n\in \Omega^{\b a}.
\]
Expand $\alpha_n$ in the basis of  $\Omega^{\b a}$, and get
\begin{equation}\label{alphav}
\alpha_v = \sum_{u\in \cal B^\lambda} c_{uv}\, z^{\deg v - \deg u-d}\, \alpha_u
\end{equation}
for some  $c_{uv}\in \k$. Now set
\[
Y(a,z)u = \sum_{v\in \cal B^{a+\lambda}} c_{uv} \,z^{\deg v - \deg u-d}\, v
\]
for any $u\in \cal B^\lambda$, and extend it by linearity to the whole $V^\lambda$.

For example, take $\lambda = 0$. Assume that $\cal B^a = \{a, Da, D^2a, \ldots\}$,  then
$\alpha_{D^m a} = \frac 1{m!}z^m \in \Omega^a_m$, and therefore
\[
Y(a, z)\1 = \sum_{m\ge 0} \frac 1{m!} \,z^m D^m a,
\]
which agrees with the identity \eqref{Y1}.

It is easy to check that~\eqref{longY} is satisf\/ied: Indeed, we have checked that it holds for
$\lambda = 0$; assuming that it holds for $\lambda = a_1 + \cdots
+ a_l$, we compute, using \eqref{alphav},
\begin{align*}
Y(a,z)Y(a_1,z_1)\cdots Y(a_l,z_l)\1 &= \sum_{u\in \cal B^\lambda} Y(a,z)u\
  \alpha_u(z_1,\ldots,z_l)\\
&= \sum_{u\in \cal B^\lambda,\,v\in \cal B^{a+\lambda}} c_{uv}z^{\deg v - \deg u-d}\,
  \alpha_u(z_1,\ldots,z_l)\, v\\
& = \sum_{v\in \cal B^{a+\lambda}} \alpha_v(z,z_1,\ldots,z_l)\, v.
\end{align*}

In order to show that the correspondence $\cal G\ni a\mapsto Y(a,z)\in
\op{Hom}(V,V\llp z\rrp)$ can be extended to a map $Y:V\to \op{Hom}(V,V\llp z\rrp)$,
we need to introduce another property of $\Omega$:

\renewcommand{\theenumi}{$\Omega$\arabic{enumi}}
\begin{enumerate}
\setcounter{enumi}{4}
\item\label{O:assoc}
If $\b a = \b a' \otimes \b a''$ for $\b a',\b a''\in T(\cal G)$, $|\b a'|=k$, $|\b a''|=l-k$,
then any function  $\alpha\in \Omega^{\b a}$ has an expansion
\begin{equation}\label{coprod}
\alpha(z_1+z,\ldots,z_k+z,z_{k+1},\ldots,z_l)=
\sum_{n\ge n_0} z^{-n-1} \sum_i \alpha'_{ni}(z_1,\ldots,z_k)\, \alpha''_{ni}(z_{k+1},\ldots,z_l)
\end{equation}
at $z = \infty$, where $\alpha'_{ni}\in \Omega^{\b a'}$ and $\alpha''_{ni}\in \Omega^{\b a''}$.
The second sum here is f\/inite.
\end{enumerate}
\renewcommand{\theenumi}{\alph{enumi}}

Note that if $k=l$, then  the expansion \eqref{coprod} just the usual Taylor
formula
\[
\alpha(z_1+z,\ldots,z_l+z) = \exp(\Delta z)\alpha(z_1,\ldots,z_l),
\]
since the left-hand side is polynomial in $z$.

\begin{lemma}\label{lem:coprod}
  Let $\Omega = \bigoplus_{\b a \in T(\cal G)} \Omega^{\b a}$ be a homogeneous space of rational
  functions, satisfying the conditions \hbox{\rm(\ref{O:1})--(\ref{O:corr})} of
  Section~{\rm \ref{sec:algtocor}}. Then it also satisfies {\rm(\ref{O:assoc})}.
\end{lemma}

Before proving this lemma, let us show how condition (\ref{O:assoc}) helps to construct the
vertex algebra structure on $V$, and hence proving Theorem~\ref{thm:V}.
Take two weights $\lambda =  a_1 + \cdots + a_l$,
$\mu = b_1 + \cdots + b_k\in \Z_+[\cal G]$,
and def\/ine the tensors  $\b a = a_1 \otimes \cdots \otimes a_l$,
$\b b = b_1\otimes\cdots\otimes b_k\in T(\cal G)$. We are going to def\/ine the action of
$Y(V^\lambda,z)$ on $V^\mu$ and then by linearity extend $Y$ to the whole~$V$.

In analogy with deriving \eqref{alphav}, we obtain from  \eqref{coprod} that every basic function $\alpha_v\in
\Omega^{\b a\otimes \b b}$ has expansion
\begin{gather}
\alpha_v(z_1+z,\ldots,z_l+z,y_1,\ldots,y_k) \nonumber\\
\qquad{}=  \sum_{u\in \cal B^\lambda,\, w\in \cal B^\mu}
c^v_{u,w}\, z^{\deg v- \deg u-\deg w}\,\alpha_u(z_1,\ldots
z_l)\, \alpha_w(y_1,\ldots,y_k),\label{expan}
\end{gather}
for some  $c^v_{u,w}\in \k$. Now we set for $u\in \cal B^\lambda$ and $w\in \cal B^\mu$
\[
Y(u,z)\,w = \sum_{v\in \cal B^{\lambda+\mu}} c^v_{u,w}\, z^{\deg v- \deg u-\deg w} \, v.
\]

To check \eqref{Yassoc}, sum \eqref{expan} over all $v\in \cal B^{\lambda+\mu}$. By \eqref{longY},
the left-hand is
\[
Y(a_1,z_1+z)\cdots Y(a_l,z_l+z)Y(b_1,y_1)\cdots Y(b_k,y_k)\1,
\]
whereas the right-hand side is, using the def\/inition of $Y(u,z)\,w$ and \eqref{longY},
\begin{gather*}
\sum_{u\in \cal B^\lambda,\, w\in \cal B^\mu} Y(u,z)w\
\alpha_u(z_1,\ldots,z_l)\,\alpha_w(y_1,\ldots,y_k) \\
\qquad{}= Y\big(Y(a_1,z_1)\cdots Y(a_l,z_l)\1, z\big)Y(b_1,y_1)\cdots Y(b_k,y_k)\1.
\end{gather*}

It remains to be seen that the map $Y:V \to \op{Hom}(V,V\llp z\rrp)$ def\/ines a structure of
vertex algebra on $V$. By the  construction, $Y$ satisf\/ies \eqref{Y1} and \eqref{DY}, and
the identity \eqref{longY} guarantees that the correlation functions for $Y$ satisfy the
\hyperlink{rat}{rationality} and \hyperlink{com}{commutativity} conditions, which, as it
was observed in Section~\ref{sec:corrfun}, are enough for $V$ to be a vertex algebra.

Note also that  $Y$ does not depend on the choice of the
bases $\cal B^\lambda$, since it depends only on the tensors $\sum_{u\in \cal B_d^\lambda}
u\otimes \alpha_u \in V_d^\lambda\otimes \Omega^{\b a}$.
\end{proof}

It is easy to see that the vertex algebra $V = V(\Omega)$ has the following universality
property:
\begin{proposition}\label{prop:uni}
Let $U$ be a vertex algebra, generated by the set $\cal G\subset U$, such that the
coalgebra of generating functions $\Omega(U)$ (given by Definition~{\rm \ref{dfn:vccf}}) is
a subspace of $\Omega$. Then there is a
unique vertex algebra homomorphism $V\to U$ that fixes $\cal G$.
\end{proposition}

\begin{remark}
In Section~\ref{sec:algtocor} we have constructed a vertex coalgebra $\Omega=\Omega(V)$ of correlation
functions of a vertex algebra $V$ (see Def\/inition~\ref{dfn:vccf}).
If we apply the construction of Theorem~\ref{thm:V} to this $\Omega$, we get $V(\Omega) =
\bigoplus_{\lambda\in \Z_+[\cal G]} V^{\lambda}$, which is the graded
deformation algebra (a.k.a.\ the Rees algebra) of~$V$.
\end{remark}

\subsection{Example: Free vertex algebra}\label{sec:free}
Clearly the conditions (\ref{O:1})--(\ref{O:corr}) are satisf\/ied for
\[
\Omega^{\b a} = \set{\alpha\in\Phi^l}{\ord_{ij}\alpha\ge -{\rm loc}(a_i,a_j) \ \forall 1\le
  i<j\le l}^{\Gamma_{\b a}}.
\]
By Proposition~\ref{prop:uni}, the resulting vertex algebra $F = F_{\rm loc}(\cal G) =
V(\Omega)$,
given by Theorem~\ref{thm:V}, has the following universal property:  any vertex algebra $U$
generated by the set $\cal G$ such that the locality of any $a,b\in \cal G$ is at
most ${\rm loc}(a,b)$, is a homomorphic image of $F$. Such a vertex algebra~$F$ is called {\em a
  free vertex algebra}. It was constructed in \cite{freecv,cfva} using dif\/ferent methods.

\subsection{Proof of Lemma~\ref{lem:coprod}}

Take some $\alpha \in \Omega^{\b a}_d$. As it is the case with any rational function with
poles at $z_i-z_j$ only,
$\alpha$ has an expansion~\eqref{coprod}. We just have to show that
$\alpha'_{ni} \in \Omega^{\b a'}$ and $\alpha''_{ni} \in \Omega^{\b a''}$.

First we show that any $\alpha''_{ni}$ in~\eqref{coprod} belongs to $\Omega^{\b a''}$. Let
$\alpha''(z_{k+1},\ldots,z_l)$ be the coef\/f\/icient of some monomial $z^{-n-1} z_1^{-n_1-1}
\cdots z_k^{-n_k-1}$ in~\eqref{coprod}. Clearly, it is enough to show that this $\alpha''\in
\Omega^{\b a''}$. The idea is that $\alpha''$ is a f\/inite linear combination of the
coef\/f\/icients of $z_1^{-m_1-1}\cdots z_k^{-m_k-1}$ in the expansion of $\alpha$ in the domain
$|z_1|> \cdots > |z_l|$, which are in $\Omega^{\b a''}$ by (\ref{O:corr1}). While this can
be shown by some manipulations with rational functions, we will use some vertex algebra
considerations.

Namely, we are going to use the free vertex algebra $F = F_{\rm loc}(\cal G)$, discussed in
Section~\ref{sec:free}. Since every function $\alpha\in\Phi^l$ satisfying (\ref{O:finite}) is a correlation function on
$F$, there is a linear functional $f: F^{\b a}_d \to \k$ such that $\alpha = \alpha_f$ is the
correlation function of $f$, given by \eqref{corrseries}.
By the \hyperlink{assoc}{associativity property} (see Section~\ref{sec:corrfun}), we have that $\alpha'' = \alpha_{f''}$ is
the correlation function of the functional $f'':F^{\b a''}_{d''}\to \k$, given by
$v \mapsto f\big((a_1(n_1)\cdots a_k(n_k)\1)(n)v\big)$, where
$d'' = d - \deg (a_1(n_1)\cdots a_k(n_k)\1)(n)$.
 Using the identity \eqref{assoc}, we see that
$(a_1(n_1)\cdots a_k(n_k)\1)(n)$ as an operator $F^{\b a''}_{d''} \to F^{\b a}_d$ can be
represented as a linear combination of words $u=a_{i_1}(m_1) \cdots a_{i_k}(m_k)\in U(F)$ for some
$m_i \in \Z$ and a permutation $\sigma = (i_1,\ldots,i_k)\in \Sigma_k$. But the
correlation function of the functional $v \mapsto f(uv)$ for such $u$ is the coef\/f\/icient of
$z_{i_1}^{-m_1-1}\cdots z_{i_k}^{-m_k-1}$ in the expansion of $\sigma\alpha$ in the domain
$|z_{i_1}|> \cdots > |z_{i_k}|> |z_{k+1}| > \cdots > |z_l|$, and therefore belongs to
$\Omega^{\b a''}$ by~(\ref{O:corr1}).

\begin{remark}
Actually, one can show that it  suf\/f\/ices  to use only words $u$ with $\sigma =1$.
\end{remark}

Now we prove that  $\alpha'_{ni} \in \Omega^{\b a'}$. Recall that  $\sigma \Omega^{\b a} =
\Omega^{\sigma \b a}$ for any permutation $\sigma \in \Sigma_l$. Apply the
permutation that reverses the order of variables to \eqref{coprod}, replace $z$ by $-z$, and
then the above argument
shows that the expansion of $\alpha(z_1,\ldots,z_k,z_{k+1}-z,\ldots,z_l-z)$ in $z$ at
$\infty$ has form
\[
\alpha(z_1,\ldots,z_k,z_{k+1}-z,\ldots,z_l-z)=
\sum_{n\ge n_0} z^{-n-1} \sum_i \~\alpha'_{ni}(z_1,\ldots,z_k)\, \~\alpha''_{ni}(z_{k+1},\ldots,z_l),
\]
where $\~\alpha'_{ni} \in \Omega^{\b a'}$. Now take another variable $w$ and
consider a f\/inite expansion
\[
\alpha(z_1 + w,\ldots,z_l+w) = \sum_j w^j \alpha^{(j)}(z_1,\ldots,z_l),
\]
where $\alpha^{(j)} =  \frac1{j!} \Delta^j\alpha\in \Omega^{\b a}$. Here we use (\ref{O:Delta})
and the fact that $\Delta$ is locally nilpotent on $\Phi^l$. Then we have
\begin{gather*}
\alpha(z_1 +w,\ldots,z_k+w, z_{k+1}+w-z,\ldots,z_l+w-z) \\
\qquad{} = \sum_j w^j \alpha^{(j)}(z_1,\ldots,z_k,z_{k+1}-z,\ldots,z_l-z)\\
\qquad{} =  \sum_{j,n,i} w^j z^{-n-1} \
\~{\big(\alpha^{(j)}\big)}'_{ni}(z_1,\ldots,z_k)\,
\~{\big(\alpha^{(j)}\big)}''_{ni}(z_{k+1},\ldots,z_l).
\end{gather*}
As we have seen,  $\big(\alpha^{(j)}\big)'_{ni}\in \Omega^{\b a'}$.
Now substitute $w = z$ in the above, and get that $\alpha'_{ni}$ is a f\/inite linear
combination of $\~{\big(\alpha^{(j)}\big)}'_{ni}$'s.
\qed

\subsection{Coalgebras of regular functions}\label{sec:caofrf}

Suppose that we are in the setup of Section~\ref{sec:cortoalg}, and that $\deg
a\ge 0 $ for any $a\in \cal G$.

\begin{theorem}\label{thm:reg0gen}
Let $\Omega_0 = \bigoplus_{\b a \in T(\cal G)} \Omega_0^{\b a}$, where
$\Omega_0^{\b a}\subset \Phi^l$, be  a homogeneous space of
functions of degree $-\sum_i \deg a_i$ for $\b a =a_1\otimes\cdots\otimes a_l\in T(\cal
G)$. Assume that for any partition $\{1,\ldots,l\}= I\sqcup J$ the component decomposition
\eqref{comp} of a function  $\alpha \in \Omega_0^{\b a}$ is
\[
\alpha = \sum_{n\ge 0} (\alpha)_n, \qquad
(\alpha)_n = \sum_j \alpha'_{-n,j} \, \alpha''_{n,j}.
\]
Assume also that
\renewcommand{\theenumi}{{\rm \roman{enumi}}}
\begin{enumerate}\itemsep=0pt
\item\label{reg0gen:1}
$\Omega_0^1=\k$;
\item\label{reg0gen:reg}
any $\alpha\in \Omega_0^{\b a}$ is $(2\deg a_1,\ldots,2\deg a_l)$-regular;
\item\label{reg0gen:ord}
$\ord_{ij}\alpha\ge -{\rm loc}(a_i,a_j)$ for every $1\le i<j\le l$;
\item\label{reg0gen:sigma}
$\sigma\Omega_0^{\b a} = \Omega_0^{\sigma\b a}$ for any permutation
$\sigma\in \Sigma_l$;
\item\label{reg0gen:nz}
for any $a_1\in\cal G$ there is a tensor $\b a = a_1 \otimes a_2\otimes \cdots \in T(\cal G)$ such that
$\Omega_0^{\b a}\neq 0$;
\item\label{reg0gen:co}
$\alpha_{0,j}'\in \Omega_0^{\b a'}$ and $\alpha_{0,j}''\in \Omega_0^{\b a''}$,
where
$\b a' = \b a(I)$ and $\b a'' = \b a(J)$.
\end{enumerate}
\renewcommand{\theenumi}{\alph{enumi}}
Let $\Omega$ be the span of all functions
$\alpha''_{n,j}$ for $n\ge 0$, so that $\alpha''_{n,j}\in\Omega^{\b a''}$.  Then $\Omega$
is a vertex coalgebra in the sense of Definition~{\rm \ref{dfn:vccf}} whose degree zero component is $\Omega_0$. The corresponding vertex
algebra $V = V(\Omega)$, given by Theorem~{\rm \ref{thm:V}}, is radical-free.
\end{theorem}

\begin{remark}
Note that the map $\Omega_0^{\b a} \to \Omega_0\otimes\Omega_0$ given by
\begin{equation}\label{co0}
 \alpha \mapsto \sum_{I\sqcup J=\{1,\ldots,l\}} \sum_j \alpha'_{0,j}\otimes \alpha''_{0,j}
\end{equation}
makes $\Omega_0$ into a coassociative cocommutative coalgebra. The dual  structure on
$V_0$ is that of an  associative commutative algebra with respect to the product $(-1)$.
\end{remark}

\begin{example}\label{ex:O}
The main example of the coalgebra $\Omega_0$ that satisf\/ies the assumptions of
Theorem~\ref{thm:reg0gen} is obtained in the following way. In the
setup of Section~\ref{sec:algtocor}, suppose that $\deg a_i\ge0$. For any $\b a = a_1\otimes \cdots\otimes
a_l\in T(\cal G)$ def\/ine
$\Omega^{\b a}_0 = \set{\alpha_f}{f:V_0\to \k, \ f(\rad V)=0}$. In particular, taking $V$
to be a free vertex algebra, introduced in Section~\ref{sec:free},  we obtain $\Omega_0^{\b a}$ being
the space of all regular $\Gamma_{\b a}$-symmetric functions $\alpha\in \Phi^l$ such that $(\alpha)_n=0$ for all $n<0$, and
$\ord_{ij}\alpha\ge -{\rm loc}(a_i,a_j)$.

Similarly, setting $\deg a = 2$ for every $a\in \cal G$ and ${\rm loc}(a_i,a_j)=4$, we can take
the space $\Omega_0^{\b a} = (R^l)^{\Gamma_{\b a}}$
of all $\Gamma_{\b a}$-invariant admissible function (see Section~\ref{sec:adm}) as another example of
a~family $\Omega_0^{\b a}$, satisfying the  assumptions of
Theorem~\ref{thm:reg0gen}.

Another similar example, that we will need in Section~\ref{sec:proofB} below, is $\Omega^{\b a}_0 =
(S^l)^{\Gamma_{\b a}}$, where $S^l\subset R^l$ is the space of admissible functions with
only simple poles.
\end{example}

\begin{proof}[Proof of Theorem~\ref{thm:reg0gen}]
Condition (\ref{O:1}) holds because of (\ref{reg0gen:nz}), and it is
easy to see that $\Omega$ satisf\/ies conditions  (\ref{O:finite}), (\ref{O:sym}) and (\ref{O:corr1}).
In order to show that $\Omega$ is indeed the vertex coalgebra generated by $\Omega_0$, we are left to check~(\ref{O:Delta}).

Take a partition $\{1,\ldots,l\} = I\sqcup J$ and set
\begin{gather*}
\Delta' = \sum_{i\in I} \partial_{z_i}, \qquad
\Delta'' = \sum_{i\in J} \partial_{z_i},\\
{\Delta^*}' = \sum_{i\in I} z^2_i\,\partial_{z_i}+(2\deg a_i)\,z_i, \qquad
{\Delta^*}'' = \sum_{i\in J} z^2_i\,\partial_{z_i}+(2\deg a_i)\,z_i.
\end{gather*}
For a function $\alpha\in \Omega_0^{\b a}$, apply $\Delta$ and $\Delta^*$ to the expansion \eqref{comp}, and get
\begin{gather*}
  0=\Delta\alpha =\sum_{n\ge 0}\sum_j \big(\Delta'\alpha_{nj}'\big)\,\alpha_{nj}''
+ \alpha_{nj}'\big(\Delta''\alpha_{nj}''\big),\\
  0=\Delta^*\alpha =\sum_{n\ge 0}\sum_j \big({\Delta^*}'\alpha_{nj}'\big)\,\alpha_{nj}''
+ \alpha_{nj}'\big({\Delta^*}''\alpha_{nj}''\big).
\end{gather*}
From this we deduce that
\begin{gather*}
  0 =(\Delta\alpha)_n = \sum_j \big(\Delta'\alpha_{nj}'\big)\,\alpha_{nj}'' +
\alpha_{n+1,j}'\big(\Delta''\alpha_{n+1,j}''\big),\\
  0 =(\Delta^*\alpha)_n = \sum_j \big({\Delta^*}'\alpha_{nj}'\big)\,\alpha_{nj}'' +
\alpha_{n-1,j}'\big({\Delta^*}''\alpha_{n-1,j}''\big),
\end{gather*}
which implies that
$\Delta'\alpha_{nj}' \in \spn\{\alpha'_{n+1,j}\}$,
$\Delta''\alpha_{nj}'' \in \spn\{\alpha''_{n-1,j}\}$,
${\Delta^*}'\alpha_{nj}' \in \spn\{\alpha'_{n-1,j}\}$, and
${\Delta^*}''\alpha_{nj}'' \in \spn\{\alpha''_{n+1,j}\}$.

Now we show that
$\rad(V)=0$. First we observe that $\rad(V)_0 = 0$, since the correlation functions of
degree 0 on $V$ being regular implies that $D^*V_1 = 0$. Now assume that there is a
homogeneous element $0\neq v\in \rad(V)$ of degree
$n>0$ and weight $\lambda = b_1 + \cdots + b_l$ for $b_i\in \cal G$. Then
there is a functional $f:V_n^\lambda\to \k$ such that $f(v)\neq 0$. Let
$\beta(z_1,\ldots,z_l) \in \Omega_n^{\b b}$ for $\b b = b_1\otimes \cdots \otimes b_l$ be
the corresponding correlation function. By the construction of $\Omega$ we can assume that $\beta$ is the coef\/f\/icient
of some monomial $w_1^{-m_1-1}\cdots w_k^{-m_k-1}$ in the power series expansion
of a function $\alpha(w_1,\ldots,w_k,z_1,\ldots,z_l)\in \Omega_0^{\b a \otimes \b b}$
in the domain $|w_1|> \cdots > |w_k|$, where $\b a = a_1 \otimes \cdots \otimes a_k$.
But then $a_1(m_1)\cdots a_k(m_k)v\neq 0$ in $V_0$, which contradicts to the fact that
$v\in \rad(V)$.
\end{proof}

\subsection{The component of degree zero}\label{sec:dg0}
Suppose $\cal G$, $T(\cal G)$,  $\rm loc$  and $\Gamma_{\b a}$ for $\b a\in T(\cal G)$ are as
in Section~\ref{sec:algtocor}. Here we prove the following fact:
\begin{theorem}\label{thm:reg0free}
Assume that for any $\b a \in T(\cal G)$ we are given a space  $\Phi^{\b a}\subset \Phi^l$, such that
the space $\Omega_0 = \bigoplus_{\b a\in T(\cal G)} \Omega_0^{\b a}$, defined by
$\Omega_0^{\b a} = (\Phi^{\b a})^{\Gamma_{\b a}}$, satisfies the assumptions of Theorem~{\rm \ref{thm:reg0gen}}. Assume also that
$\Phi^{\b a}\Phi^{\b b}\subset \Phi^{\b a\otimes \b b}$ for any $\b a, \b b\in T(\cal G)$.
Let $V = V(\Omega)$ be the vertex algebra constructed in Theorem~{\rm \ref{thm:reg0gen}}. Then
$V_0$ is isomorphic to a polynomial algebra.
\end{theorem}
Note that the spaces $\Omega_0$ given in Example~\ref{ex:O} are all obtained in this
way.

Before proving this theorem, we need to establish certain property of the algebra $V_0$.
We know that $V_0$ is an associative commutative algebra, graded by weights:
$V_0 = \bigoplus_{\lambda\in \Z_+[\cal G]} V_0^{\lambda}$.
 Let $X = \bigoplus_{\lambda\neq 0} V_0^\lambda$ be the augmentation
ideal in $V_0$. Consider the symmetrized tensor product $\op{Sym}^2_X X = \big(X \otimes_X
X\big)_{\Sigma_2}$. There is the canonical homomorphism $\mu:\op{Sym}^2_X X\to X^2$
def\/ined by $\mu(x\otimes y) = xy$.

\begin{lemma}\label{lem:V0}
The map $\mu:\op{Sym}^2_X X\to X^2$ is an isomorphism.
\end{lemma}

\begin{proof}
Clearly, $\mu$ is surjective. To prove that it is also injective,
suppose that $\sum_i u_iv_i = 0$ in $X^2\subset V_0$ for some homogeneous $u_i,v_i \in
X$. We need to show that
$\sum_i u_i\otimes v_i = 0$ in  $\op{Sym}^2_X X$.

The tensor product $\op{Sym}^2_X X$ is graded by $\Z_+[\cal G]$.
Therefore, it is enough to check that $f\big(\sum_i u_i\otimes v_i\big)=0$ for any
homogeneous linear functional
$f: \op{Sym}^2_X X\to \k$. Assume that $\wt f = \lambda = a_1 + \cdots + a_l\in \Z_+[\cal G]$.
For a
non-trivial partition $P = \{P_1,P_2\} \in \bar{\cal P}_2$, set $\lambda' = \sum_{i\in
  P_1}a_i$ and $\lambda'' =
\sum_{i\in P_2}a_i$. Then $f$ can be pulled back to a functional on
$V_0^{\lambda'}\otimes V_0^{\lambda''}$. Since both $V_0^{\lambda'}$ and $V_0^{\lambda''}$
are f\/inite-dimensional, we can write this functional as
$\sum_j f'_j \otimes f''_j$ for some $f'_j:V_0^{\lambda'}\to \k$ and $f''_j:V_0^{\lambda''}\to \k$.

Set $\b a = a_1 \otimes \cdots\otimes a_l$ and $\b a' = \b a(P_1)$, $\b a'' = \b a(P_2)$ as
in Section~\ref{sec:algtocor}.
Let $\alpha'_j\in \Omega_0^{\b a'}$ and $\alpha''_j\in \Omega_0^{\b a''}$
be the correlation functions of $f'_j$ and $f''_j$ respectively. Set
\[
\alpha(P) = \sum_j \alpha'_j\, \alpha''_j.
\]

Denote $\Gamma = \Gamma_{\b a}$. We claim that the
functions $\alpha(P)\in \Phi^{\b a}$ for $P\in \bar{\cal P}_2$ satisfy the properties of
Proposition~\ref{prop:comp0} and also
$\alpha(\sigma P) = \alpha(P)$ for any $\sigma\in \Gamma$.

Note that one of the assumptions of Theorem~\ref{thm:reg0gen} was that $(\alpha')_d = (\alpha'')_d =
0$ for \mbox{$d<0$}, therefore $(\alpha(Q))_0(P)$ is the leading term in the expansion~\eqref{comp}
of a function~$\alpha(Q)$.
The condition~\eqref{alphaPQ} follows from the fact that $ab\otimes cd=ac\otimes bd$ in $\op{Sym}^2_X
X$ for every  $a,b,c,d\in X$.

So by Proposition~\ref{prop:comp0} there exists a function $\alpha\in \Phi^{\b a}$ such that
$(\alpha)_0(P) = \alpha(P)$ for any partition $P \in \bar{\cal P}$.
Replacing $\alpha$ by $|\Gamma|\inv\sum_{\sigma\in \Gamma}\sigma\alpha$ we can assume that
$\alpha\in \Omega_0^{\b a}$.
Then $\alpha$ is a~correlation function of a
linear functional $h:V_0^\lambda\to \k$, such that $h(u_i v_i) = f(u_i\otimes v_i)$ for any pair $u_i$,
$v_i$. Therefore, $f\big(\sum_i u_i\otimes v_i\big) = h\big(\sum_i u_i v_i\big) = 0$.
\end{proof}

\begin{proof}[Proof of Theorem~\ref{thm:reg0free}]
Recall that the augmentation ideal $X = \bigoplus_{l>0}X_l$ of $V_0$ is graded, where
\[
X_l = \bigoplus_{\b a\in T(\cal G),\, |\b a|=l} V_0^{\b a}.
\]
 For $v\in X_l$ we will call $l=|v|$ the length of~$v$.
Choose a homogeneous basis  $\cal X \subset X$ of $X$ modulo~$X^2$. We want to show that $V_0 \cong \k[\cal X]$.
Note that we can extend the grading on $\cal X$ to the grading on~$\k[\cal X]$.

Consider the canonical map $\f:\k[\cal X]\to V_0$, that maps every element $x\in\cal X$
into itself.  We need to show that $\f$ is an isomorphism. It is easy to see that $\f$ is
surjective~-- this follows from the fact that for  f\/ixed length $l$, we have $X_l\cap X^k
= 0$ for $k\gg 0$.

Let $\bar X = \cal X\k[\cal X]$ be the augmentation ideal of $\k[\cal X]$.
Consider the map $\psi:\bar X^2\to \op{Sym}^2_XX$ that maps a monomial
$x_1\cdots x_k$ to $x_1 \otimes \f(x_2 \cdots x_k)$ for $x_i\in \cal X$.
Note that the space $\op{Sym}^2_XX$  is graded by the length.

The restriction  $\f:\bar X_i\to X_i$ is an isomorphism for the minimal $i$, because then
$X_i = \bar X_i = \spn\set{x\in \cal X}{|x|=i}$.
Assume we have established that $\f:\bar X_i\to X_i$ is an isomorphism for $i\le l-1$.
Then  $\psi:\bar X^2_i\to \big( \op{Sym}^2_XX \big)_i$ is an isomorphism for $i\le l$.
Combining this with the isomorphism of Lemma~\ref{lem:V0}, we get an isomorphism
$\bar X^2_i \cong X^2_i$ for $i\le l$. But if $p \in \bar X_l$ is such that $\f(p)=0$,
then $p\in \bar X^2$, since $\cal X\cup \{\1\}$ is linearly independent modulo~$X^2$,
therefore we must have $p \equiv 0$ and $\bar X_l \cong X_l$.
\end{proof}

\section{OZ vertex algebras}\label{sec:OZ}
\subsection{Some notations}\label{sec:rho}
Assume we have a vertex algebra $V$ graded as $V =V_0\oplus V_2\oplus V_3\oplus \cdots$.
First of all recall (see Section~\ref{sec:def}) that $V_0$ is an associative commutative algebra under the
operation $ab = a(-1)b$ and  $V$ is a vertex algebra over
$V_0$. Indeed, $V$ is a $V_0$-module under the action $av=a(-1)v$ for
$a\in V_0$ and $v\in V$, and the identity (\hyperlink{iv}{V4}) of
Def\/inition~\ref{dfn:vertex} implies that this action commutes with the vertex algebra
structure on $V$.
The component $V_2$ is a
commutative (but not associative in general) algebra with respect to the product $ab=a(1)b$,
equipped with an invariant symmetric bilinear form $\form ab = a(3)b$.

Let $A\subset V_2$ be a subspace such that $A(3)A\subseteq \k\1\subseteq V_0$ and
$A(1)A\subseteq A$. Set as before $T(A) = \set{a_1 \otimes \cdots \otimes a_l\in
A^{\otimes l}}{a_1,\ldots,a_l\in A}$.  Denote by $V'$ the graded dual space of $V$.

For a
tensor $\b a = a_1 \otimes \cdots \otimes a_l\in T(A)$ and a linear
functional $f\in V'$ let $\alpha= \alpha_f(z_1,\ldots,z_l)$ the correlation function, given by
\eqref{corrseries}. Then $\alpha\in R^l$. Indeed, $\alpha$ is regular, since $D^*V_1=0$ (see
Section~\ref{sec:regular}), has $(\alpha)_n=0$ for $n<0$ or $n=1$ because $V_n=0$ for these $n$ (see
Section~\ref{sec:components}), has poles of order at most 4 since ${\rm loc}(a,b)=4$ for any $a,b\in A$ and
has $\rho^{(k)}_{ij}\in R^{l+k/2}$ for $k=-2,-4$
by the \hyperlink{assoc}{associativity} property of
correlation functions. Note that $\rho_{ij}^{(-4)}\alpha$ does not depend on~$z_i$,~$z_j$,
since $a_i(3)a_j \in \k\1$. This def\/ines a map
\begin{equation}\label{phi}
\phi:V'\otimes A^{\otimes l} \to R^l,
\end{equation}
such that $\phi(f,\b a) = \sigma \phi(f,\sigma\b a)$ for any $\b a \in A^{\otimes l}$,
 $f\in V'$ and $\sigma\in \Sigma_l$.

 There is an obvious action of the symmetric groups and of
$\k^\times$ on $T(A)$. Set $S(A) = T(A)_{\Sigma\times
\k^\times} = PT(A)_\Sigma$. Denote by
$\Omega^{\b a} = \set{\phi(f,\b a)}{f\in V'}\subset R^l$ the space of all
correlation functions corresponding to $\b a$.

Recall that for $\b a = a_1\otimes \cdots\otimes a_l\in T(A)$ we have considered the space
\[
V^{\b a} = \spn_\k\bigset{a_1(m_1)\cdots a_l(m_l)\1}{m_i\in \Z},
\]
 so that $V^{\b a} \cong \big(\Omega^{\b a}\big)'$.
If $\b a = \b b$ in $S(A)$, then $V^{\b a} = V^{\b b}$. Set also
\begin{equation}\label{filtrV}
V^{(l)}  = \spn_\k \bigset{a_1(m_1)\cdots a_k(m_k)\1}{a_i\in A, m_i\in
  \Z, k\le l},
\end{equation}
so that we have a f\/iltration $\k\1 = V^{(0)} \subseteq
V^{(1)}\subseteq V^{(2)}\subseteq \cdots\subseteq V$.
Denote $V^{(l)}_d = V^{(l)}\cap V_d$.

Let $\cal G\subset A$ be a linear basis of $A$. Then $\Z_+[\cal G]$ can be identif\/ied with a subset of $S(A)$ by $a_1 + \cdots + a_l =
a_1 \otimes \cdots\otimes a_l$ for $a_i \in \cal G$. Every tensor $\b a \in T(A)$ of length $l$ can be expanded as
$\b a = \sum_i k_i \b g_i$ for $k_i\in \k$ and $\b g_i \in T(\cal G)$. This implies that
\begin{equation}\label{VcalG}
V^{(l)} = \bigcup_{\lambda \in \Z_+[\cal G], \, |\lambda|\le l} V^{\lambda}.
\end{equation}

Next we def\/ine the maps $r_{ij}^{(1)}, r_{ij}^{(3)}:T(A)\to T(A)$ by
\begin{gather*}
r_{ij}^{(1)}a_1\otimes\cdots \otimes a_l  =
a_1\otimes\cdots \otimes a_{i-1}\otimes a_ia_j\otimes \cdots
\otimes{\^a}_j\otimes\cdots a_l,\\
r_{ij}^{(3)}a_1\otimes\cdots \otimes a_l  = \form{a_i}{a_j}\
a_1\otimes\cdots \otimes a_{i-1}\otimes \^a_i\otimes \cdots
\otimes\^a_j\otimes\cdots a_l.
\end{gather*}

It follows from  the \hyperlink{assoc}{associativity property} of Section~\ref{sec:corrfun} that
\begin{equation}\label{rhophi}
\rho_{ij}^{(-2)} \phi(f,\b a) = \phi\big(f,r_{ij}^{(1)}\b a\big),\qquad
\rho_{ij}^{(-4)} \phi(f,\b a) = \phi\big(f,r_{ij}^{(3)}\b a\big)
\end{equation}
for all $f\in V'$.

 Def\/ine a partial ordering on $T(A)$ by writing $\b b \prec \b a$ if
 $\b b = r_{ij}^{(k)}\b a$, and taking the transitive
 closure.

\subsection{Virasoro element}\label{sec:vir}
Now we investigate what happens when an element $\omega\in A\subset
V_2$ is a Virasoro element of $V$ (see Section~\ref{sec:def}).
Recall that in this case
$\omega(0)a = Da$ and $\omega(1)a = (\deg a)\, a$ for every homogeneous $a\in
V$, see \eqref{virconf}, and therefore  $\frac12 \omega$ is an idempotent in the Griess algebra $V_2$.

Consider a tensor  $\b a = a_1\otimes\cdots\otimes a_{l-1}\otimes \omega \in T(A)$, and
set
\begin{equation}\label{b}
\b b = r_{il}^{(1)}\b a = 2\, a_1\otimes\cdots\otimes a_{l-1}, \ \
\b b_i = r_{il}^{(3)}\b a =
\form \omega{a_i}\, a_1\otimes\cdots\otimes \^a_i\otimes\cdots\otimes
a_{l-1}.
\end{equation}
Let $f\in V'_0$ be a linear functional. Denote $\alpha(z_1, \ldots, z_l) =
\phi(f, \b a)$,  $\beta(z_1,\ldots,z_{l-1})=\phi(f, \b b)$ and $\beta_i(z_1, \ldots,\^z_i, \ldots,z_{l-1}) =
\phi(f, \b b_i)$.

Def\/ine an operator $\cal E:\Phi^{l-1}\to \Phi^{l-1}$ by
\[
\cal E = \sum_{i=1}^{l-1}-z_i\inv\partial_{z_i}+2z_i^{-2},
\]
and let the shift operator $T:\Phi^{l-1}\to \Phi^l$ be given by
$T(f(z_1,\ldots,z_{l-1})) = f(z_1-z_l,\ldots,z_{l-1}-z_l)$.

\begin{proposition}\label{prop:Evir}
Suppose that $V$ is generated by $A\subset V_2$ as a vertex
algebra. Then
an element $\omega\in A$ is a Virasoro element of $V$ if and only if
\begin{equation}\label{Evir}
\alpha = \frac 12 T\cal E\beta + \sum_{i=1}^{l-1} (z_i-z_l)^{-4}
\beta_i
\end{equation}
for every  $\b a = a_1\otimes\cdots\otimes a_{l-1}\otimes \omega \in
T(A)$ and $f\in V'_0$.
\end{proposition}

\begin{remark}
If $\form \omega\omega \neq 0$, we can choose a basis $\omega\in \cal G\subset A$ so
that $\omega$ is orthogonal to the rest of basic elements. Then it is
enough to check~\eqref{Evir} only for $a_1\otimes\cdots\otimes a_{l-1}\in
T(\cal G)$, and the second sum runs over the indices $i$ such that $a_i=\omega$.
\end{remark}

\begin{proof}
If $\omega \in V_2$ is a Virasoro element, then by \eqref{com} we have
$\ad{a(m)}{\omega(-1)} = (m+1)\,a(m-2) +
\delta_{m,3}\form \omega a$ for any $a\in A$, which implies
\[
\ad{Y(a,z)}{\omega(-1)} = (2z^{-2}-z\inv\partial_z)\,Y(a,z) + \form \omega a \, z^{-4}\1.
\]
Therefore,
\begin{gather*}
\alpha(z_1,\ldots, z_{l-1},0) =
f\big(Y(a_1,z_1)\cdots Y(a_{l-1},z_{l-1})\, \omega(-1)\1\big)\\
\phantom{\alpha(z_1,\ldots, z_{l-1},0)}{} = \cal Ef\big(Y(a_1, z_1)\cdots Y(a_{l-1}, z_{l-1})\1\big) \\
\phantom{\alpha(z_1,\ldots, z_{l-1},0)=}{} + \sum_{i=1}^{l-1}  \form \omega{a_i}\, z_i^{-4}
f\big(Y(a_1,z_1)\cdots \^{Y(a_i,z_i)}\cdots Y(a_{l-1},z_{l-1})\, \1\big)\\
\phantom{\alpha(z_1,\ldots, z_{l-1},0)}{}=\frac 12 \cal E
\beta(z_1,\ldots,z_{l-1}) +  \sum_{i=1}^{l-1}
z_i^{-4} \beta_i(z_1, \ldots,\^z_i, \ldots,z_{l-1}),
\end{gather*}
and we get  \eqref{Evir} since
$\alpha(z_1,\ldots,z_l) = \alpha(z_1-z_l,\ldots,z_{l-1}-z_l,0)$ by
Proposition~\ref{prop:regular}.

Conversely, in order to see that $\omega \in A$ is a Virasoro element,
we need to show that $\op{ad} \omega(1):U(V)\to U(V)$ is the grading
derivation and $\op{ad} \omega(0):U(V)\to U(V)$ coincides with
$D$. Since~$A$ generates $V$ as a vertex algebra, the operators
$a(n)$ for  $a\in A$,  $n\in \Z$, generate $U(V)$ as an associative algebra, and
therefore it is enough to verify commutation relations between
$\omega(m)$ and $a(n)$ for $m=0,1$ and $n\in \Z$. Using \eqref{com} this
amounts to checking the identities
\[
\omega(0) a = Da \qquad \text{and}  \qquad \omega(1)a = 2a
\]
for any $a\in A$. Note that we also have $\omega(2) a = 0$ since
$V_2=0$  and
$\omega(3)a = \form \omega a\,\1$ since $\omega\in A$.
Using~\eqref{qs}, these identities are equivalent to
\begin{equation}\label{vir2check}
a(0)\omega = Da,\qquad a(1)\omega = 2a,\qquad \forall \, a\in A.
\end{equation}
Setting $a_{l-1}=a$ and $z_{l-1}=z$ we expand
\[
\alpha(z_1,\ldots,z_{l-1},0) = \sum_{n\ge 0}
z^{-4+n}\alpha_n(z_1,\ldots,z_{l-2}),
\]
where $\alpha_n=f\big(Y(a_1,z_1)\cdots Y(a_{l-2},z_{l-2})\, a(3-n)\omega\big)$.
It follows that \eqref{vir2check} is equivalent to
\begin{gather*}
\alpha_2(z_1,\ldots,z_{l-2}) = f\big(Y(a_1,z_1)\cdots
Y(a_{l-2},z_{l-2})\, a(1)\omega\big) \\
\phantom{\alpha_2(z_1,\ldots,z_{l-2})}{} = 2 f\big(Y(a_1,z_1)\cdots Y(a_{l-2},z_{l-2})\, a\big) =
\beta(z_1,\ldots,z_{l-2},0),\\
 \alpha_3(z_1,\ldots,z_{l-2}) = f\big(Y(a_1,z_1)\cdots
 Y(a_{l-2},z_{l-2})\, a(0)\omega\big) \\
\phantom{\alpha_3(z_1,\ldots,z_{l-2})}{}= f\big(Y(a_1,z_1)\cdots Y(a_{l-2},z_{l-2})\, Da\big) = \frac12 \frac{\partial \beta}{\partial z_{l-1}}\Big|_{z_{l-1}=0},
\end{gather*}
which easy follows from~\eqref{Evir}.
\end{proof}

\subsection[The operator $\cal E$]{The operator $\boldsymbol{\cal E}$}
In this section we show that the operator $\cal E$ preserves the property of being admissible (see
Section~\ref{sec:adm}).

\begin{proposition}\label{prop:Eadm}
For an admissible function $\beta \in R^{l-1}$,   $l\ge3$, set  $\alpha =T\cal E\beta$.
\begin{enumerate}\itemsep=0pt
\item\label{Eadm:il}
For any $1\le i<l$ we have
\[
\rho^{(-1)}_{il}\alpha = - \partial_{z_i}\beta,\qquad
\rho^{(-2)}_{il}\alpha = 2\beta
\]
and $\rho^{(k)}_{il}\alpha  = 0$ for $k<-2$.

\item\label{Eadm:ij}
For any $1\le i<j<l$ we have
\[
\rho^{(-2)}_{ij}\alpha = T\cal E\rho^{(-2)}_{ij}\beta + 2(z_j-z_l)^{-4} \rho_{ij}^{(-4)}\beta,\qquad
\rho^{(-4)}_{ij}\alpha = T\cal E\rho^{(-4)}_{ij}\beta
\]
and $\rho^{(k)}_{ij}\alpha  = 0 \quad \text{for} \ \ k<-4$.

\item\label{Eadm:reg}
$\alpha \in R^l$.
\end{enumerate}
\end{proposition}

\begin{proof}(\ref{Eadm:il})
 Since $\beta$ and
$\partial_{z_i}\beta$ are translation-invariant, we have
\[
T\cal E\beta =
\sum_{i=1}^{l-1} \Big(-(z_i-z_l)\inv\partial_{z_i} + 2(z_i-z_l)^{-2}\Big)\beta(z_1,\ldots,z_{l-1}),
\]
therefore
\[
\alpha = 2 (z_i-z_l)^{-2}\beta - (z_i-z_l)\inv
\partial_{z_i}\beta + O\big((z_i-z_l)^0\big).
\]

\noindent
(\ref{Eadm:ij})
Assume for simplicity that $i=1$ and $j=2$. Expand
\[
\beta = (z_1-z_2)^{-4}\beta_{-4}(z_3,\ldots,z_{l-1}) + (z_1-z_2)^{-2}\beta_{-2}(z_2,\ldots,z_{l-1})+\cdots.
\]
Using that $\ad{\cal E}{z_1-z_2} = z_1\inv z_2\inv (z_1-z_2)$, we get
\[
\cal E (z_1-z_2)^k
= (z_1-z_2)^k(\cal E+k\, z_1\inv z_2\inv).
\]
So we compute
\begin{align*}
\cal E(z_1-z_2)^{-2}\beta_{-2}
&= (z_1-z_2)^{-2}\Big(\cal E' + 2 z_1^{-2} - 2 z_1\inv z_2\inv\Big)\beta_{-2}\\
&=   (z_1-z_2)^{-2} \cal E'\beta_{-2} + O\big((z_1-z_2)^{-1}\big),\\
\cal E(z_1-z_2)^{-4}\beta_{-4}
&= (z_1-z_2)^{-4}\Big(\cal E'' + 2 z_1^{-2} +2 z_2^{-2} - 4 z_1\inv z_2\inv\Big)\beta_{-4}\\
&= (z_1-z_2)^{-4}\cal E''\beta_{-4} + 2 (z_1-z_2)^{-2} z_2^{-4}\beta_{-4} +
O\big((z_1-z_2)^{-1}\big),
\end{align*}
where $\cal E' = \sum_{i=2}^{l-1} z_1\inv \partial_{z_1}- 2z_1^{-2}$, \
$\cal E'' =  \sum_{i=3}^{l-1} z_1\inv \partial_{z_1}- 2z_1^{-2}$.

\medskip\noindent
(\ref{Eadm:reg})
First we show that $\alpha$ is regular.
Set $\Delta_1^* = \sum_{i=1}^{l-1} z_i^2\partial_{z_i}+4z_i$ and
$\Delta^*=\Delta_1^*+z_l^2\partial_{z_l}+4z_l$.  It is enough to check that
$\Delta_1^*(\cal E\beta)(z_1,\ldots,z_{l-1}) = 0$. Indeed, in this case
set $w_i = z_i -z_l$ for $1\le i\le l-1$, and get
\begin{align*}
\Delta^* T(\cal E\beta) &=
\(\Big(\sum_{i=1}^{l-1} z_i^2\partial_{w_i} + 4 z_i\Big) -
\Big(\sum_{i=1}^{l-1}z_l^2\partial_{w_i}\Big) + 4z_l\)\cal E\beta(w_1, \ldots,w_{l-1})\\
&= \(\Delta_1^* + z_l\Big(4l+2\sum_{i=1}^{l-1}w_i\partial_{w_i}\Big) \)\cal E\beta(w_1, \ldots,w_{l-1})\\
&=(4l + 2\deg \cal E\beta) \, z_l\, \cal E\beta(w_1, \ldots,w_{l-1}) = 0,
\end{align*}
since $\deg \cal E\beta = \deg \beta - 2 = -2l$.
So we compute $\ad{z^2\partial_z+4z}{-z^{-2}\partial_z+2z} = 3\partial_z$, hence
$\ad{\Delta_1^*}{\cal E} = 3\sum_{i=1}^{l-1}\partial z_i$,
and therefore, using Proposition~\ref{prop:regular}, we get $\Delta_1^*\cal E \beta = \cal E\Delta_1^* \beta + 3 \sum_i \partial_{z_i}\beta =
0$.

In order to f\/inish the proof of  (\ref{Eadm:reg}) we need only
to show that for every partition
$\{1,\ldots,l\} = I\sqcup J$ the expansion \eqref{comp} of $\alpha$ has form
\[
\alpha =(\alpha)_0 + \sum_{n\ge2}(\alpha)_n
\]
so that $(\alpha)_0 =
\sum_j\alpha'_{0,j}\, \alpha''_{0,j}$ for $\alpha'_{0,j}\in R^{|I|}$ and $\alpha''_{0,j}\in
R^{|J|}$.

We prove this statement by induction on $l$. If $l=3$, then
$\beta=k(z_1-z_2)^{-4}$ for
$k\in \k$ and then $\alpha = 2k\, (z_1-z_2)^{-2}(z_1-z_3)^{-2}(z_2-z_3)^{-2} \in R^3$ and the
statements (\ref{Eadm:il}) and (\ref{Eadm:ij}) are obviously true. So assume
that $l\ge 4$.

Without loss of generality we can assume that $l\in J$.
Write the expansion \eqref{comp} for $\beta$ as $\beta = (\beta)_0 + \sum_{n\ge2}(\beta)_n$
where $(\beta)_n = \sum_j\beta'_{n,j}\, \beta''_{n,j}$.
Note that both $(\beta)_0$ and $\sum_{n\ge2}(\beta)_n$ are admissible.
Then
\[
\alpha =  \sum_{n\ge0,\,n\neq 1}T\cal E (\beta)_n
=\sum_{n,j} \big(T\cal E'\beta'_{n,j}\big)\beta''_{n,j}
+\beta'_{n,j}\big(T\cal E''\beta''_{n,j}\big),
\]
where  $\cal E' = \sum_{i\in I} -z_i\inv\partial_{z_i}+2z_i^{-2}$ and
$\cal E'' = \sum_{i\in J\backslash\{l\}} -z_i\inv\partial_{z_i}+2z_i^{-2}$.
By induction, $T\cal E'\beta'_{0,j}\in R^{|I|+1}$ and
$T\cal E''\beta''_{0,j}\in R^{|I|}$, therefore $T\cal E(\beta)_0 \in R^l$.

We are left to show that $\big(T\cal E(\beta)_n\big)_m=0$ for $n\ge2$ and $m\le1$.
Observe that $T\cal E'\beta'_{n,j}$ does not have
pole at $z_l$, therefore
\[
\big(\big(T\cal E'\beta'_{n,j}\big)\beta''_{n,j}\big)_m = 0
\]
for $m<n$, and the claim follows.
\end{proof}

\subsection[Explicit formulae for $\Omega_0^{a}$ for small $a$]{Explicit formulae for $\boldsymbol{\Omega_0^{\b a}}$ for small
$\boldsymbol{\b a}$}\label{sec:smalla}

Let $\b a = a_1\otimes\cdots\otimes a_l\in T(A)$. It follows from \eqref{rhophi} that
if $f:V_0^{\b a}\to \k$ is such that $f(V^{(l-1)})=0$, then the corresponding correlation
function $\phi(f,\b a)\in \Omega_0^{\b a}$ has only simple poles. The smallest such
function is
\[
\prod_{1\le i<j\le 5} (z_i-z_j)\inv \in R^5,
\]
therefore for $l\le 4$ the space of correlation functions $\Omega_0^{\b a}$ has dimension
1. We have $\Omega_0^1=\k$,   $\Omega_0^a = 0$, and for $l=2,3,4$, $\Omega_0^{\b a} =
\k\alpha$, where $\alpha$ is as follows:
\begin{align*}
l=2: \ \  &\alpha= \form{a_1}{a_2}\, (z_1-z_2)^{-4}\\
l=3: \ \ &\alpha = \form{a_1}{a_2a_3}\, (z_1-z_2)^{-2}(z_1-z_3)^{-2}(z_2-z_3)^{-2}\\
l=4: \ \
&\alpha = \form{a_1}{a_2}\form{a_3}{a_4}\, (z_1-z_2)^{-4}(z_3-z_4)^{-4}\\
&\!+ \form{a_1}{a_3}\form{a_2}{a_4}\, (z_1-z_3)^{-4}(z_2-z_4)^{-4}\\
&\!+ \form{a_1}{a_4}\form{a_2}{a_3}\, (z_1-z_4)^{-4}(z_2-z_3)^{-4}\\
&\!+ \form{a_1a_2}{a_3a_4}\, (z_1-z_2)^{-2}(z_3-z_4)^{-2}(z_1-z_3)\inv(z_1-z_4)\inv(z_2-z_3)\inv(z_2-z_4)\inv\\
&\!+ \form{a_1a_3}{a_2a_4}\, (z_1-z_3)^{-2}(z_2-z_4)^{-2}(z_1-z_2)\inv(z_1-z_4)\inv(z_2-z_3)\inv(z_3-z_4)\inv\\
&\!+ \form{a_1a_4}{a_2a_3}\, (z_1-z_4)^{-2}(z_2-z_3)^{-2}(z_1-z_2)\inv(z_1-z_3)\inv(z_2-z_4)\inv(z_3-z_4)\inv.\!
\end{align*}

\section[The algebra $B$]{The algebra $\boldsymbol{B}$}\label{sec:B}
Now let $A$ be a commutative algebra with a symmetric invariant bilinear form
$\formdd$. Denote by $\op{Aut}A$ its the group of automorphisms, i.e.\ the linear
maps that preserve the product and the form on $A$. In this section we prove the following
theorem:

\begin{theorem}\label{thm:B}
There exists a vertex algebra
$B = B_0 \oplus B_2 \oplus B_3 \oplus \cdots$, generated by
$A\subset B_2$, so  that
\begin{enumerate}\itemsep=0pt
\item\label{B:A}
$a(1)b = ab$ and $a(3)b = \form ab$
for any $a,b\in A$;
\item\label{B:Vir}
if $\frac12\omega\in A$ is a unit of $A$, then $\omega$ is a Virasoro element of $B$;
\item\label{B:canon}
$\op{Aut}A\subset \op{Aut}B$.
\item\label{B:B0free}
If $\dim A = 1$ or, if $A$ has a unit $\frac12 \omega$, $\dim A/\k \omega = 1$, then $B_0 =\k$; otherwise,
$B_0$ is  isomorphic to the polynomial algebra in infinitely many
  variables.
\end{enumerate}
\end{theorem}

\begin{remark}
In fact one can show that the vertex algebra $B$ can be obtained
as $B = \^B/\rad \^B$, where $\^B$ is the vertex algebra generated by the space $A$
subject to  relations (\ref{B:A}) of Theorem~\ref{thm:B} and condition $\^B_1=0$.
\end{remark}

Before constructing the algebra $B$ and proving Theorem~\ref{thm:B}, let us show how it implies the
main result of this paper.

\subsection{Proof of Theorem~\ref{main}}\label{sec:proof}
Assume that the form $\formdd$ on $A$ is non-degenerate. Take an arbitrary algebra homomorphism
$\chi: B_0 \to \k$. If a f\/inite group of automorphisms $G\subset \op{Aut}A$ was specif\/ied,
choose $\chi$ to be $G$-symmetric. (Here we use that $\op{char}\k=0$.)
By  Proposition~\ref{prop:forms}, $\chi$ def\/ines a $\k$-valued
symmetric bilinear form $\formdd_\chi$ on $B$, such that ${\form ab}_\chi = \chi\big(a(-1)^*b\big)$.
It is easy to see that the form $\formdd_\chi$ coincides with the form $\formdd$ on $A$. Indeed,
for $a,b\in A$ we have, using that $\chi(\1) = 1$,
\[
\form ab_\chi = \chi\big(a(-1)^*b\big) = a(3)b \, \chi(\1) = \form ab.
\]

Since $\chi$ is multiplicative, we have $\ker_\chi \subset \rad\formdd_\chi$.
Now we set
\[
V = B/\rad \formdd_\chi.
\]
Since $V_0 = \k\1$, Proposition~\ref{prop:rad}\ref{rad:ideals} implies that $V$ is simple.
This proves Theorem~\ref{main}\ref{main:emb}, whereas (\ref{main:vir}) and (\ref{main:aut})
of Theorem~\ref{main}  follow from  (\ref{B:Vir}) and (\ref{B:canon}) of Theorem~\ref{thm:B}.

\subsection[Constructing $B_0$]{Constructing $\boldsymbol{B_0}$}\label{sec:B0}

Let us f\/ix a linear basis $\cal G$ of $A$, such that if $A$ has a unit
$\frac 12 \omega$, then $\omega\in \cal G$.
First we construct the spaces
\begin{equation}\label{Bfiltr}
\k\1=B_0^{(0)}\subseteq B_0^{(1)} \subseteq \cdots \subseteq B_0,
\end{equation}
def\/ined  by \eqref{filtrV}.
Recall that $T(A) = \set{a_1\otimes a_2\otimes \cdots }{a_i\in A}\subset \bigoplus_{l \ge
  1} A^{\otimes l}$ and $S(A) = T(A)_{\Sigma \times \k^\times}$.
For any $\b a \in T(A)$ of length $|\b a|\le l$
we will construct a subspace $B_0^{\b a}\subset B_0^{(l)}$ and the dual space of admissible
correlation functions $\big(B_0^{\b a}\big)^*=\Omega_0^{\b a}\subset R^l$, so that the following properties will hold:

\renewcommand{\theenumi}{B\arabic{enumi}}
\begin{enumerate}\itemsep=0pt

\item\label{B:SA}$B_0^{\b a} =B_0^{\b b}$ if $\b a = \b b$ in $S(A)$.
\item\label{B:prec} $B_0^{\b b} \subseteq B_0^{\b a}$ whenever $\b b \prec \b a$.
\item\label{B:phi}  There is a map
\[
\phi:\big(B_0^{\b a}\big)^*\otimes A^{\otimes l} \to \Omega_0^{\b a},
\]
for $\b a \in T(A)$,  $|\b a|=l$, satisfying  \eqref{rhophi} and
$\phi(f,  \b a) = \sigma\phi(f,\sigma\b a)$
for any permutation $\sigma\in\Sigma_l$.

\item\label{B:vir}
If $\frac12 \omega\in A$ is a unit, then for any $\b b\in T(A)$, $|\b b|=l-1$, $\b a = \b b \otimes
\omega$, \, $\b b_i = r^{(3)}_{ij}\b a$ and $f:B^{(l)}\to \k$,
\[
\phi(f, \b a) = T\cal E \phi(f,\b b) + \sum_{i=1}^{l-1} (z_i-z_l)^{-4} \phi(f, \b b_i).
\]

\item\label{B:gr}
$\displaystyle{B_0^{(l)}/B_0^{(l-1)} = \bigoplus_{\lambda \in \Z_+[\cal G\backslash \{\omega\}]}
B_0^\lambda / B_0^{(l-1)}}$.

\end{enumerate}
\renewcommand{\theenumi}{\alph{enumi}}
In addition we want $\Omega_0^{\b a}$ to satisfy the conditions (\ref{reg0gen:1})--(\ref{reg0gen:co}) of Theorem~\ref{thm:reg0gen}.

The map $\phi$ in (\ref{B:phi}) is going to be the same as in \eqref{phi}.
As in Section~\ref{sec:algtocor}, the condition~(\ref{B:phi}) implies that $\Omega^{\b a} \subset
R^{\Gamma_{\b a}}$. The condition (\ref{B:vir}) is the same as in Proposition~\ref{prop:Evir}.
The property~(\ref{B:SA}) justif\/ies the
 notation $B_0^\lambda = B_0^{\b a}$ used in  (\ref{B:gr}), whenever $\lambda = g_1 +
 \cdots + g_l\in \Z_+[\cal G]$ and $\b a = g_1\otimes \cdots \otimes g_l\in T(\cal G)$.
The property (\ref{B:gr}) is a special case of~\eqref{VcalG}.

We are constructing $B_0^{(l)}$ by induction on $l$, starting from
$B_0^{(0)}=\k\1$ and $\Omega_0^1 = \k$. Assume that $B_0^{(m)}$, $\Omega_0^{\b a}$ and
$\phi:\big(B_0^{(m)}\big)^*\otimes A^{\otimes m}\to R^m$ are already constructed  for
$m=|\b a| \le l-1$.

\subsubsection*{Constructing $\boldsymbol{\Omega_0^{\b g}}$}
Take a basic tensor $\b g = g_1\otimes \cdots\otimes g_l\in T(\cal G)$.
We def\/ine the space $\Omega_0^{\b g}\subset R^l$ in the
following way:
If $g_i \neq \omega$ for all $1\le i\le l$, then set
\[
\Omega_0^{\b g} = \Biggset{\alpha\in (R^l)^{\Gamma_{\b g}}}{
\begin{gathered}
\exists f:B^{(l-1)}_0\to \k \ \ s.\,t. \ \  \rho^{(-k-1)}_{ij}\alpha = \phi\big(f,r^{(k)}_{ij}\b g\big)\\
\forall \ 1\le i<j\le l, \ k = 1\ \text{or} \ 3
\end{gathered}
}.
\]
It is not clear a priori why $\Omega_0^{\b g}\neq 0$. This is a part of the statement of
Proposition~\ref{prop:Omega} below. Note also that by f\/ixing some $0\neq \alpha \in \Omega_0^{\b g}$, the space
$\Omega_0^{\b g}$ can be described as the space of functions that dif\/fer from $\alpha$ by
an admissible $\Gamma_{\b g}$-symmetric function with only simple poles.

If $g_l =\omega$, then set
$\b b = r_{il}^{(1)}\b g, \, \b b_i = r_{il}^{(3)}\b g\in T(\cal G)$
as in \eqref{b},
and then def\/ine $\Omega_0^{\b g}$ to be set of all functions $\alpha\in R^l$ given by
\eqref{Evir}, where $\beta = \phi^{l-1}(f, \b b)$, $\beta_i = \phi^{l-2}(f, \b b_i)$ for all
 $f:B_0^{(l-1)}\to \k$. Note that $\alpha\in R^l$ due to Proposition~\ref{prop:Eadm}\ref{Eadm:reg}.

Finally, if $g_i = \omega$ for some $1\le i\le l-1$, then set
\[
\Omega_0^{\b g} =  (i\,l)\Omega_0^{(i\,l)\b g},
\]
for the  transposition  $(i\, l)\in \Sigma_l$.

It is immediately clear that  for any permutation $\sigma \in \Sigma_l$ we have
\begin{equation}\label{symm}
\Omega_0^{\b g} = \sigma \Omega_0^{\sigma \b g}.
\end{equation}
The following proposition gives another crucial property of the functions $\Omega_0^{\b g}$.

\begin{proposition}\label{prop:Omega}
For any linear functional $f:B_0^{(l-1)}\to \k$ and any tensor $\b g =g_1\otimes\cdots\otimes
g_l\in T(\cal G)$ there is a function $\alpha \in \Omega_0^{\b g}$ such that
$\rho_{ij}^{(-k-1)} \alpha = \phi(f, r_{ij}^{(k)} \b g)$.
\end{proposition}

\begin{proof}
In the case when $\b g$ does not contain $\omega$, Proposition~\ref{prop:ML} guarantees that there exists
an admissible  function $\alpha\in R^l$ such that $\rho_{ij}^{(-k-1)} \alpha = \phi(f,
r_{ij}^{(k)} \b g)$, since the functions $\alpha_{ij}^{(k)} = \phi(f,r_{ij}^{(-k-1)}\b g)$
obviously satisfy the condition \eqref{rhoalpha}. Now take $|\Gamma_{\b g}|\inv
\sum_{\sigma\in \Gamma_{\b g}} \sigma\alpha\in \Omega_0^{\b g}$.

Now suppose  that $\b g$ contains $\omega$. Using \eqref{symm}, we can assume without
loss of generality, that $\b g = a_1\otimes \cdots\otimes a_{l-1}\otimes \omega$.
Let $\b b = r_{il}^{(1)}\b g$ and $\b b_s = r_{sl}^{(3)}\b g$ for $1\le s< l$ as in
\eqref{b}, and  set
$\beta = \phi(f,\b b)$, \,  $\beta_s= \phi(f,\b b_s)$. Then $\alpha=\phi(f,\b g)\in
\Omega^{\b g}_0$ is def\/ined by
\eqref{Evir}. By Proposition~\ref{prop:Eadm}\ref{Eadm:il}, we get $\rho_{il}^{(-2)}\alpha = \beta$ and
$\rho_{il}^{(-4)}\alpha =\beta_i$ for all $1\le i<l$.

Now consider the case when  $1\le i<j<l$.
Applying $\rho^{(-4)}_{ij}$ to \eqref{Evir} and using \break
Proposition~\ref{prop:Eadm}\ref{Eadm:ij}, we get
\[
\rho^{(-4)}_{ij}\alpha = T\cal E\rho^{(-4)}_{ij}\beta + \sum_{s\neq i,j} (z_s-z_l)^{-4}\,
\rho^{(-4)}_{ij}\beta_s = \phi(f,r_{ij}^{(3)}\b g).
\]
To do the same with $\rho_{ij}^{(-2)}$ we notice that
\[
r^{(3)}_{jl} r^{(1)}_{ij}\b g = r^{(3)}_{ij}\b b.
\]
Indeed,  assuming that $i,j=1,2$ to simplify notations, we get
\begin{gather*}
r^{(3)}_{2l} r^{(1)}_{12}\b g  =
r^{(3)}_{2l} (a_1a_2)\otimes a_3\otimes \cdots\otimes a_{l-1}\otimes \omega\\
\phantom{r^{(3)}_{2l} r^{(1)}_{12}\b g}{}
= \form{\omega}{a_1a_2}\, a_3\otimes \cdots\otimes a_{l-1}
= 2\form{a_1}{a_2}\, a_3\otimes \cdots\otimes a_{l-1} = r_{12}^{(3)}\b b.
\end{gather*}
Now we apply $\rho^{(-2)}_{ij}$ to \eqref{Evir} and compute, using Proposition~\ref{prop:Eadm}\ref{Eadm:ij},
\begin{equation*}
\rho^{(-2)}_{ij}\alpha = T\cal E\rho^{(-2)}_{ij}\beta +
(z_j-z_l)^{-4}\rho_{ij}^{(-4)}\beta + \sum_{s\neq i,j} (z_s-z_l)^{-4}\,
\rho^{(-2)}_{ij}\beta_s = \phi(f,r_{ij}^{(1)}\b g).
 \tag*{\qed}
\end{equation*}
  \renewcommand{\qed}{}
\end{proof}

\subsubsection*{Constructing $\boldsymbol{B_0^{(l)}}$}
For $\b g \in T(\cal G)$ set $B_0^{\b g} = \big(\Omega_0^{\b g}\big)^*$. For a permutation
$\sigma \in \Sigma_l$ we identify $B_0^{\b g}$ with $B_0^{\sigma\b g}$ by setting $(b,
\alpha) = (b,\sigma\inv\alpha)$ for $b\in B_0^{\sigma\b g}$ and $\alpha \in \Omega_0^{\b g}$,
since $\sigma\inv\alpha \in \Omega^{\sigma\b g}$ by~\eqref{symm}. Thus for a weight $\lambda =
g_1+\cdots + g_l \in \Z_+[\cal G]$ we can denote $B_0^\lambda= B_0^{\b g}$, where $\b g =
g_1 \otimes \cdots \otimes g_l\in T(\cal G)$.

By the construction, for any function $\alpha \in \Omega_0^{\b g}$ there is a linear functional
$f:B_0^{(l-1)}\to \k$ such that $\rho^{(m)}_{ij}\alpha = \phi\big(f, r^{(-m-1)}_{ij}\b
g\big)$ for all $1\le i<j\le l$ and $m=-2,-4$. Take a tensor $\b g \succ \b a\in T(A)$. There is a
map $\rho_{\b g\b a}: \Omega_0^{\b g} \to \Omega_0^{\b a}$, which is an iteration of the maps
$\rho^{(m)}_{ij}$, so that
\begin{equation}\label{rhoga}
\rho_{\b g\b a}\alpha = \phi(f, \b a).
\end{equation}
We note that the restriction of $f$ on $B_0^{\b a}$ is uniquely def\/ined by $\alpha$, so
$\rho_{\b g \b a}$ is well def\/ined by \eqref{rhoga}.

By Proposition~\ref{prop:Omega}, the map $\rho_{\b g\b a}:\Omega_0^{\b g}\to
\Omega_0^{\b a}$ is surjective, therefore we have an embedding \break
$\rho_{\b g\b a}^*: B_0^{\b a}\hookrightarrow B_0^{\b g}$.
Set
\begin{equation}\label{B0l}
B_0^{(l)} = \bigg(B_0^{(l-1)}\oplus\bigoplus_{\substack{\lambda\in \Z_+[\cal G]\\|\lambda|=l}} B_0^\lambda
\bigg)\bigg/
\spn\Bigset{a-\rho_{\b g\b a}^*a}{a\in B_0^{\b a},\, T(\cal G) \ni \b g \succ \b a\in T(A)}.
\end{equation}
In other words, we identify the space $B_0^{\b a}\subset B_0^{(l-1)}$ with
the subspace $\rho_{\b g\b a}^*(B_0^{\b a})\subset B_0^{\b g}$ for
$\b g \succ \b a \in T(A)$. So we have $B_0^{(l-1)}\subset B_0^{(l)}$
and $B_0^\lambda \subset B_0^{(l)}$ for any $\lambda \in \Z_+[\cal G]$ of
length $l$.

For $\b g \in T(\cal G)$, $|\b g|=l$, we def\/ine the map
$\phi(\,\cdot\,, \b g):\big(B^{(l)}_0\big)^*\to \Omega_0^{\b g}\subset R^l$
in the following way. Since $\big(B^{\b g}_0\big)^*\cong \Omega_0^{\b g}$, the restriction of a functional
$f:B^{(l)}_0\to \k$ to $B^{\b g}_0$ can be identif\/ied with a function in $\Omega_0^{\b g}$, and we set
$\phi(f, \b g) = f\big|_{B_0^{\b g}}$. Then we extend $\phi$ by linearity to the map
$\phi:\big(B^{(l)}_0\big)^*\otimes A^{\otimes l}\to R^l$.

\subsubsection*{Verif\/ication of  (\ref{B:SA})--(\ref{B:gr})}
The properties (\ref{B:SA}) and (\ref{B:phi})
are clear. For a tensor $\b a \in T(A)$, $|\b a|=l$, consider
the restriction $\phi(\,\cdot\,,\b a):\big(B^{(l)}_0\big)^*\to
R^l$. Denote its image by $\Omega_0^{\b a}\subset R^l$. The the dual map
$\phi(\,\cdot\,,\b a)^*:B_0^{\b a}\to B_0^{(l)}$ is an embedding of
the dual space $B_0^{\b a} = \big(\Omega_0^{\b a}\big)^*$ into
$B_0^{(l)}$. This establishes~(\ref{B:prec}).

If $\b a\in T(A)$ ends by $\omega$, and  $\b b$ and $\b b_i$ are def\/ined as in
\eqref{b}, then any
function $\alpha \in \Omega^{\b a}$ is uniquely def\/ined by its coef\/f\/icients
$\rho_{il}^{(-2)}\alpha = \beta\in\Omega^{\b b}$
and $\rho_{il}^{(-4)}\alpha = \beta_i\in \Omega^{\b b_i}$ by the formula~\eqref{Evir}.
This implies (\ref{B:vir}).
Note also that in this case
\[
B^{\b a} \subseteq B^{\b b} + B^{\b b_1} + \cdots + B^{\b b_{l-1}} \subset B^{(l-1)}.
\]
The condition (\ref{B:gr}) follows from \eqref{B0l}.

\subsubsection*{Verif\/ication of conditions (\ref{reg0gen:1})--(\ref{reg0gen:co}) of
  Theorem~\ref{thm:reg0gen}}

The only conditions of Theorem~\ref{thm:reg0gen} that require verif\/ication are (\ref{reg0gen:nz}) and
(\ref{reg0gen:co}). If we assume that the form $\formdd$ on $A$ is non-degenerate, then
for any $a\in \cal G$ there is $b\in \cal G$ such that
$\form ab\neq 0$, and then $\Omega_0^{a\otimes b} = \k\, (z_1-z_2)^{-4}\neq0$, which
proves (\ref{reg0gen:nz}).
Another argument, that does not use non-degeneracy of $\formdd$, can be found in the proof of
Theorem~\ref{thm:B}\ref{B:B0free} in Section~\ref{sec:proofB} below.

To prove (\ref{reg0gen:co}),
consider a partition $\{1,\ldots,l\}= I\sqcup J$ and set $\b a' = \b
a(I), \, \b a'' = \b a(J) \in T(\cal G)$ (in the notations of Section~\ref{sec:algtocor}). Since $\alpha\in R^l$, we have
$(\alpha)_n =\sum_j \alpha'_{-n,j}\,\alpha''_{n,j} =0$ for $n<0$ or $n=1$ by Def\/inition~\ref{dfn:adm},
and
$\alpha'_{0,j},\alpha''_{0,j}$ are admissible by Proposition~\ref{prop:adm}\ref{adm:poly}.
Denote $\b b = r^{(k)}_{ij} \b a$ and $\b b' = \b b(I)$. To show that $\alpha'_{0,j}\in \Omega_0^{\b a'}$ we need to show that
$\rho^{(-k-1)}_{ij}\alpha'_{0,j}\in \Omega_0^{\b b'}$,
for any $i,j\in I$ and $k=1$ or 3.  But this follows from (\ref{reg0gen:co}) applied to
$\rho^{(-k-1)}_{ij}\alpha$, since  we obviously have
\[
\big(\rho^{(-k-1)}_{ij}\alpha\big)_0 = \sum_s \big(\rho^{(-k-1)}_{ij}\alpha'_{0,j}\big)(\alpha''_{0,j}).
\]
Similarly, we check that $\alpha''_{0,j}\in \Omega_0^{\b a''}$.

\subsection{Proof of Theorem~\ref{thm:B}}\label{sec:proofB}

We apply Theorem~\ref{thm:reg0gen} to the space of functions $\Omega_0$
constructed in Section~\ref{sec:B0}, and obtain a~vertex algebra
\[
\^B = \bigoplus_{n\ge0, \, n\neq 1} \^B_n, \qquad
\^B_n =\bigoplus_{\lambda \in \Z_+[\cal G]}\^B^\lambda_n.
\]
Note that we have $\^B_1=0$ due to the fact that $\Omega_0$ consists
of admissible functions which
do not have components of degree 1, see Def\/inition~\ref{dfn:adm}.

From the construction of $B_0$ we see that $\^B_0^\lambda =
B_0^\lambda$ for any $\lambda \in \Z_+[\cal G]$. Recall that $\^B_0$ is the associative commutative
 algebra,  dual to  the coalgebra $\Omega_0$ with respect to the coproduct~\eqref{co0}. But then $B_0$ is also an associative commutative algebra,
since $B_0^\lambda \subset B_0$, and we have a surjective algebra homomorphism $\pi: \^B_0
\to B_0$. Its kernel $\ker \pi$ is an ideal in $\^B_0$, which  by Proposition~\ref{prop:rad}\ref{rad:ideals}
can be extended to an ideal $\ol{\ker \pi} \subset \^B$. So we f\/inally set
\[
B = \^B / \ol{\ker \pi}.
\]

The condition (\ref{B:A}) holds by the construction: indeed, given $\b a = a_1\otimes
a_2\otimes a_3\otimes \cdots \in T(A)$, and a correlation function $\alpha\in \Omega^{\b
  a}$, the coef\/f\/icient $\rho_{12}^{(-2)}\alpha$ is the correlation function corresponding to
$(a_1a_2)\otimes a_3\otimes \cdots$, but by the \hyperlink{assoc}{associativity} condition
of Section~\ref{sec:corrfun} it must be the correlation function for $\big(a_1(1)a_2\big)\otimes
a_3\otimes \cdots$, which implies that $a_1(1)a_2 = a_1a_2$. The equality $a_1(3)a_2 =
\form{a_1}{a_2}$ is established in the same way.

If $\frac12\omega\in A$ is a unit, then any correlation function $\alpha \in \Omega_0^{\b
   a\otimes \omega}$ is given by the formula \eqref{Evir}, therefore, $\omega$ is a Virasoro
   element by  Proposition~\ref{prop:Evir},  thus proving (\ref{B:Vir}).
Note also that the construction of $B$ was canonical, which establishes
(\ref{B:canon}).

\subsubsection*{Proof of Theorem~\ref{thm:B}\ref{B:B0free}}

We need to introduce another vertex algebra.
Let $F = V(\Omega(F))$ be the vertex algebra obtained by the construction of
Theorem~\ref{thm:reg0gen} from the space of functions
\[
\Omega(F)_0 = \bigoplus_{\b g \in T(\cal G\backslash \{\omega\})} \Omega(F)_0^{\b g},
\]
where $\Omega(F)_0^{\b g} = (S^l)^{\Gamma_{\b g}}$ is the
space of $\Gamma_{\b g}$-symmetric admissible functions with only simple poles
(see Example~\ref{ex:O}).
By Theorem~\ref{thm:reg0free} the algebra $F_0$ is polynomial.
Note that the algebra $F$ is what Theorem~\ref{thm:B} would yield  if instead of $A$ one would take the
space  $\spn\big\{\cal G\backslash \{\omega\}\big\}$
with zero product and form.

On the other hand, the vertex algebra $B$ inherits a
f\/iltration \eqref{Bfiltr} from $B_0$. Consider
the associated graded algebra $\op{gr}B = \bigoplus_{l\ge 1}
B^{(l)}/B^{(l-1)}$. This is indeed a vertex algebra, since all vertex
algebra identities (see Def\/inition~\ref{dfn:vertex}) are homogeneous.
By (\ref{B:gr}) we have
\[
\op{gr}B = \bigoplus_{\lambda \in \Z_+[\cal G \backslash\{\omega\}]}(\op{gr}B)^\lambda.
\]
From the construction in Section~\ref{sec:B0} it follows that the coalgebra of correlation functions
$\Omega(\op{gr}B)_0$ is the same as $\Omega(F)_0$, therefore $\Omega(\op{gr}B) =
\Omega(F)$ and hence  by Proposition~\ref{prop:uni} there is a vertex algebra
isomorphism $\eta:(\op{gr}B) \to F$, which yields algebra isomorphism
$\eta_0:(\op{gr}B)_0 = \op{gr}(B_0)\to F_0$.
But since a polynomial algebra cannot have
non-trivial deformations, we must have $B_0 \cong F_0$ as associative commutative algebras.

We are left with estimating the size of $F_0$. Take some
$\b g = a_1\otimes
\cdots\otimes a_l\in T(\cal G\backslash \{\omega\})$ and let $\Gamma = \Gamma_{\b g}$. A
function  $\alpha\in \Omega_0^{\b g}=(S^l)^\Gamma$ can have a pole at $z_i-z_j$ only if
$a_i\neq a_j$. Therefore, since $\deg \alpha = -2l$, if  $\cal G$ has no more than one
element other than $\omega$,  then $(S^l)^\Gamma = 0$ for $l>0$, and hence
  $F_0 = \k$.

Now assume that $\cal G$ has at least two elements other then $\omega$, say $a$ and $b$.
Denote by $S_0^l\subset S^l$  the space of indecomposable admissible functions
with only simple poles. Then the span of the generators of degree $\b g$ of $F_0$ is
 isomorphic to $(S^l_0)^\Gamma$.  We claim that for $l$ large enough there is $\b g\in T(\cal
G\backslash \{\omega\})$,  $|\b g|=l$, such that the  $(S^l_0)^\Gamma\neq 0$.
This would imply that $F_0$ is a polynomial algebra in inf\/initely many variables.

Indeed, there are
inf\/initely many bipartite 4-regular connected graphs that remain connected after a
removal of any two edges. Let $G$ be such a graph with vertices
$u_1,\ldots,u_k,v_{k+1},\ldots,v_l$, so that an edge can only connect some $u_i$ with some
$v_j$. The incidence matrix of this graph is an  $l\times l$ symmetric regular matrix
$\sf S=\{s_{ij}\}_{i,j=1}^l$ (see Example~\ref{ex:S}), def\/ined
so that $s_{ij}=-1$ whenever $G$ has an edge connecting $u_i$ and $v_j$ for some $1\le i\le
k< j\le l$, the rest
of the entries being 0. Then
\[
0\neq |\Gamma|\inv\sum_{\sigma\in\Gamma}\sigma\pi(\sf S) \in (S^l_0)^\Gamma=\Omega_0^{\b g},
\]
where $\pi(\sf S)$ is as in \eqref{pi} and $\b g =  a\otimes \cdots \otimes a
\otimes b\otimes\cdots\otimes b$.

\pdfbookmark[1]{References}{ref}
\LastPageEnding

\end{document}